\documentclass[12pt,a4paper,reqno]{amsart}
\usepackage[T1]{fontenc}
\usepackage{bm}
\usepackage{bbm}
\usepackage{amssymb}
\usepackage[abbrev]{amsrefs}
\usepackage[dvipsnames] {xcolor}
\usepackage{hyperref}
\usepackage[shortlabels] {enumitem}

\AtBeginDocument{\def\MR#1{}}
\newtheorem{thm}{}[section]
\theoremstyle{plain}
\newtheorem{theorem}[thm]{Theorem}
\newtheorem{corollary}[thm]{Corollary}
\newtheorem{lemma}[thm]{Lemma}
\newtheorem{proposition}[thm]{Proposition}
\theoremstyle{definition}
\newtheorem{definition}[thm]{Definition}
\newtheorem{example}[thm]{Example}
\theoremstyle{remark}
\newtheorem{remark}[thm]{Remark}

\numberwithin{equation}{section}
\allowdisplaybreaks
\DeclareMathOperator{\spn}{span}
\DeclareMathOperator{\supp}{supp}
\DeclareMathOperator{\sgn}{sign}
\newcommand{\udf}{{\bm{\varphi}}}
\newcommand{\isep}{\mathrel{{.}\,{.}}\nobreak}
\newcommand{\littletaller}{\mathchoice{\vphantom{\big|}}{}{}{}}
\newcommand\restr[2]{{\left.\kern-\nulldelimiterspace #1 \littletaller \right|_{#2}}}
\newcommand{\abs}[1]{\left\lvert#1\right\rvert}
\newcommand{\norm}[1]{\left\lVert#1\right\rVert}
\newcommand{\enbrace}[1]{\left\lbrace#1\right\rbrace}
\newcommand{\enpar}[1]{\left(#1\right)}
\newcommand{\Ft}{\mathcal{F}}
\newcommand{\Gt}{\mathcal{G}}
\newcommand{\It}{\mathcal{I}}
\newcommand{\Et}{\mathcal{E}}
\newcommand{\Ot}{\mathcal{O}}
\newcommand{\St}{\mathcal{S}}
\newcommand{\Ut}{\mathcal{U}}
\newcommand{\XB}{\mathcal{X}}
\newcommand{\YB}{\mathcal{Y}}
\newcommand{\BB}{\mathbb{B}}
\newcommand{\FF}{\mathbb{F}}
\newcommand{\NN}{\mathbb{N}}
\newcommand{\Ss}{\mathbb{S}}
\newcommand{\XX}{\mathbb{X}}
\newcommand{\YY}{\mathbb{Y}}
\newcommand{\ZZ}{\mathbb{Z}}
\newcommand{\Id}{\mathrm{Id}}
\newcommand{\Ind}{\mathbbm{1}}
\newcommand{\ee}{\bm{e}}
\newcommand{\ww}{\bm{w}}
\newcommand{\tq}{\colon}
\newcommand{\coo}{c_{00}}
\newcommand{\co}{c_0}
\newcommand{\uu}{{\bm{u}}}
\newcommand{\ZB}{{\mathcal{Z}}}
\newcommand{\cG}{{\mathcal{G}}}
\newcommand{\cM}{{\mathcal{M}}}

\begin{document}
\title[The role of symmetry in the TGA]{The role of symmetry in the performance of the Thresholding Greedy Algorithm}
\author[J. L. Ansorena]{Jos\'e L. Ansorena}
\address{Department of Mathematics and Computer Sciences\\
Universidad de La Rioja\\
Logro\~no\\
26004 Spain}
\email{joseluis.ansorena@unirioja.es}
\author[M. Berasategui]{Miguel Berasategui}
\address{Miguel Berasategui\\
UBA - Pab I, Facultad de Ciencias Exactas y Naturales\\
Universidad de Buenos Aires\\
(1428), Buenos Aires, Argentina}
\email{mberasategui@dm.uba.ar}
\author[A. Marcos]{Alejandro Marcos}
\address{Department of Mathematics and Computer Sciences\\
Universidad de La Rioja\\
Logro\~no\\
26004 Spain}
\email{alejandro.marcos@unirioja.es}
\thanks{J.\@ L.\@ Ansorena acknowledges the support of the Spanish Ministry for Science, Innovation, and Universities under Grant PGC2018-095366-B-I00 for \emph{Functional Analysis Techniques in Approximation Theory and Applications (TAFPAA)}.}
\thanks{A. Marcos acknowledges support from the University of La Rioja under an FPI predoctoral contract associated with the above-mentioned grant PGC2018-095366-B-I00.}
\begin{abstract}
This paper continues the study of conditional symmetric systems initiated in \cite{AnsorenaMarcos2026} and solves several problems posed there. While the interaction between unconditionality and symmetry is well known in the Banach setting following the work in \cites{AABCO2024, AnsorenaMarcos2026}, we broaden the scope and address similar questions for spreading systems in the more general setting of quasi-Banach spaces. Our aim is twofold. On the one hand, we characterize which spreading systems of Banach spaces are unconditional. On the other hand, we provide examples of quasi-Banach spaces equipped with conditional systems with symmetry-related properties to illustrate the connections between symmetry, greedy-like properties, and the Schauder basis property in the absence of local convexity.
\end{abstract}
\maketitle
\section{Introduction}\noindent
The theory of symmetric systems, that is, minimal systems equivalent to all their permutations, has been mainly studied under the assumption that the systems are unconditional. However, instances of conditional symmetric systems appear naturally in the literature (see, e.g., \cite{AABCO2024}*{Example 3.9}). The authors of \cite{AABCO2024} proved that the unconditionality of symmetric systems is closely related to the behavior of the Thresholding Greedy Algorithm (TGA for short) relative to the system. In fact, combining their advances with those achieved by the authors of \cite{AnsorenaMarcos2026} yields the following result.

\begin{theorem}[\cite{AnsorenaMarcos2026}*{Theorem 3.6}]\label{thm:AABCOM}
Let $\XB$ be a fundamental minimal system of a Banach space. If $\XB$ is symmetric, then the following are equivalent.
\begin{enumerate}[label=(\roman*)]
\item $\XB$ is unconditional.
\item $\XB$ is unconditional for constant coefficients.
\item $\XB$ is bidemocratic.
\item $\XB$ is superdemocratic.
\item\label{it:Mk} $\XB$ is a Markushevich basis.
\item\label{it:Sc} $\XB$ is a Schauder basis.
\end{enumerate}
\end{theorem}

If a minimal system is symmetric, then it is spreading, that is, equivalent to all its subsequences (see \cite{Singer1962} for a proof in the unconditional setting and \cite{AnsorenaMarcos2026} for a proof where unconditionality is not assumed). The link between spreadingness and symmetry is so strong that, in many cases, the only property we require of unconditional symmetric systems is their subsymmetry, that is, unconditionality and spreadingness. In fact, it was believed that subsymmetry forces symmetry until Garling \cite{Garling1968} disproved it. Given this connection, it is natural to ask for an analogue of Theorem~\ref{thm:AABCOM} that does not require symmetry and assumes spreadingness instead. In this paper, we tackle this task and give a complete answer. Before stating this result, we point out that the summing basis of $\co$ is a conditional spreading Schauder basis. Thus, items \ref{it:Mk} and \ref{it:Sc} of Theorem~\ref{thm:AABCOM} cannot be included in the corresponding characterization of unconditional spreading systems.

\begin{theorem}\label{thm:dem-spread-equivalences}
Let $\XB$ be a fundamental minimal system of a Banach space $\XX$. If $\XB$ is spreading, then the following are equivalent.
\begin{enumerate}[label=(\roman*)]
\item\label{cor:dem-spread-equivalences:unc} $\XB$ is unconditional.
\item\label{cor:dem-spread-equivalences:bidem} $\XB$ is bidemocratic.
\item\label{cor:dem-spread-equivalences:superdem} $\XB$ is superdemocratic.
\item\label{cor:dem-spread-equivalences:ucc} $\XB$ is unconditional for constant coefficients.
\end{enumerate}
\end{theorem}

Since every Markushevich basis has a Schauder subsequence (see \cite{KadPel1965}*{Theorem 3}), being a Schauder basis and being a Markushevich basis are equivalent properties in the locally convex and spreading context (see \cite{AnsorenaMarcos2026}*{Corollary 2.10}). We complement Theorem~\ref{thm:dem-spread-equivalences} by describing those spreading systems that fail to be Markushevich bases. Specifically, we will prove the following.

\begin{proposition}\label{prop:spread_not_markushevich}
Let $\XB=(x_n)_{n=1}^\infty$ be a spreading fundamental minimal system of a Banach space $\XX$. If $\XB$ is not a Markushevich basis, then $\XB$ converges weakly to some $y_0\in\XX\setminus\enbrace{0}$. Moreover, if we set $y_n=x_n-y_0$ for all $n\in\NN$, then the sequence $(y_n)_{n=1}^\infty$ is a weakly null unconditional basis of its closed linear span, and $(y_n)_{n=0}^\infty$ is a subsymmetric basis of $\XX$.
\end{proposition}

Another line of research starting from Theorem~\ref{thm:AABCOM} involves determining the role of local convexity within it. Every characterization of unconditionality displayed in Theorem~\ref{thm:AABCOM} and Theorem~\ref{thm:dem-spread-equivalences} can be directly translated to the more general setting of quasi-Banach spaces. So, it is only natural to fix two of these features and study whether, assuming either symmetry or spreadingness and dropping local convexity, one implies the other. An important counterexample that answers some of these questions negatively appeared in \cite{AABB2026}*{Section 4}. Specifically, the quasi-Banach space constructed there has a symmetric Markushevich basis that fails to be a Schauder basis under any rearrangement, let alone an unconditional basis. Moreover, this Markushevich basis is almost greedy and fails to be bidemocratic. Thus, symmetric fundamental systems that fail to be unconditional but still perform well relative to the TGA do exist. In this paper, we construct a different example with even better behavior.

\begin{theorem}\label{thm:BeraA}
There is a quasi-Banach space $\XX$ with a symmetric Markushevich basis $\XB$ that is bidemocratic and quasi-greedy but fails to be a Schauder basis under any rearrangement.
\end{theorem}

We also study the existence of systems that fail to be Markushevich bases, but are symmetric or, at least, spreading, and have suitable greedy-like properties. Since any quasi-greedy basis is a Markushevich basis (see \cite{AABW2021}*{Corollary 4.5}), to obtain results in this direction, we need to lower our expectations of how the TGA performs relative to the system. Our contribution here is the following.

\begin{theorem}\label{thm:BeraB}
There is a quasi-Banach space with a symmetric fundamental minimal system that is superdemocratic, semi-greedy, and truncation quasi-greedy but fails to be a Markushevich basis.
\end{theorem}

The paper is organized as follows. In Section~\ref{sec:Background}, we recall the necessary background on minimal systems and greedy approximations. In Section~\ref{sec: preliminary}, we record some auxiliary results that will be used later on. Section~\ref{sec:spreading-quasi-Banach} studies spreading fundamental systems in quasi-Banach spaces and contains the proof of Theorem~\ref{thm:BeraB}. In Section~\ref{sec:spreading-Banach}, we turn to the locally convex case. We prove Theorem~\ref{thm:dem-spread-equivalences} and obtain further results on spreading fundamental systems. In Section~\ref{sec:sym-qBanach}, we develop counterexamples that yield Theorem~\ref{thm:BeraA} and obtain additional results on symmetric systems.
\section{Terminology and background}\label{sec:Background}\noindent
The symbol $\delta_{j,k}$ denotes the Kronecker delta associated with the indices $j$ and $k$. Given $a\in\ZZ$, we define the set
\[
\ZZ_{\ge a}=\enbrace{n\in\ZZ \tq n\ge a}.
\]
For $j$, $k\in\ZZ$, $[j \isep k]$ will denote the integer interval with endpoints $j$ and $k$. The symbol $[A]_{<\infty}$ stands for the set consisting of all finite subsets of $A$. If $A$ is infinite, then $[A]_{<\infty}$ endowed with the inclusion order is a directed set. Given nonempty sets $A$, $B\in[\ZZ]_{<\infty}$, we say that $A<B$ if $\max(A)<\min(B)$.

We will denote by $\coo(J)$ the vector space of all eventually null scalar families modeled over a set $J$. We say that a family $(b_j)_{j\in J}$ in $[0,\infty]$ \emph{dominates} another family $(a_j)_{j\in J}$ in $[0,\infty]$, denoted by $a_j\lesssim b_j$ for $j\in J$, if there is a constant $C > 0$ such that $a_j\le C b_j$ for all $j\in J$. If both $a_j\lesssim b_j$ and $b_j\lesssim a_j$ for $j\in J$, the families are said to be \emph{equivalent}, and we put $a_j\approx b_j$ for $j\in J$. Given two nonnegative sequences $(a_n)_{n=1}^\infty$ and $(b_n)_{n=1}^\infty$ we say that the series $\sum_{n=1}^\infty a_n$ and $\sum_{n=1}^\infty b_n$ are equivalent, and we put
\[
\sum_{n=1}^\infty a_n \sim \sum_{n=1}^\infty b_n
\]
if their sequences of partial sums are equivalent.

We write $\Pi(A)$ for the group of permutations of a set $A$. Similarly, if $A\subset\ZZ$, $\Ot(A)$ will denote the set of all increasing maps from $A$ into itself. Note that a function $\psi\in\Ot(A)$ is fully determined by its range. Finally, we denote by $\It(A)$ the set of all one-to-one maps from $A$ into itself.

We refer the reader to \cite{KPR1984} for the basics of locally bounded metrizable topological vector spaces when local convexity is not assumed. A \emph{quasi-norm} on a vector space $\XX$ over the real or complex field $\FF$ is a map $\norm{\cdot}\colon \XX \to [0, \infty)$ such that
\begin{enumerate}[label=(\roman*)]
\item\label{QN:a} $\norm{ x}=0$ only when $x=0$,
\item\label{QN:b} $\norm{t \, f}=\abs{t} \norm {f}$ for all $t\in\FF$ and all $f\in\XX$, and
\item there is a constant $\kappa\ge 1$ so that for all $f$ and $g\in \XX$ we have
\begin{equation*}
\norm{f+g} \le \kappa\enpar{\norm{f} +\norm{g}}.
\end{equation*}
\end{enumerate}
The optimal constant $\kappa$ will be called the \emph{modulus of concavity of $\XX$}.

Given $0<p\le 1$, a $p$-norm is a function satisfying \ref{QN:a}, \ref{QN:b}, and
\begin{enumerate}[label=(\roman*),resume]
\item $\norm{f+g}^p \le \norm{f}^p +\norm{g}^p$ for all $f$ and $g\in \XX$.
\end{enumerate}
Any $p$-norm is a quasi-norm with modulus of concavity at most $2^{1/p-1}$. Conversely, by the Aoki--Rolewicz theorem, any quasi-norm with modulus of concavity at most $2^{1/p-1}$ is equivalent to a $p$-norm (see, e.g., \cite{KPR1984}).

A quasi-norm $\norm{\cdot}$ on $\XX$ induces a linear topology on the space. If $\XX$ endowed with this topology is a complete topological vector space, we call $\XX$ a \emph{quasi-Banach space}. A $p$-Banach space (Banach space if $p=1$) is a quasi-Banach space equipped with a $p$-norm. The symbol $\Ss_\XX$ stands for the unit sphere of $\XX$. Given a family $\XB=(x_j)_{j\in J}$ in $\XX$, $[\XB]$ and $[x_j \tq j\in J]$ both denote its closed linear span. The family $\XB$ in $\XX$ dominates another family $\YB=(y_j)_{j\in J}$ in a quasi-Banach space $\YY$ if the scalar family
\[
\norm{\sum_{j\in J} a_n \, x_n}, \quad (a_j)_{j\in J}\in\coo(J)
\]
dominates the corresponding family associated with $\YB$. Equivalently, $\XB$ dominates $\YB$ if there exists a bounded linear map $T\colon\XX\to\YY$ such that $T(x_j)=y_j$ for every $j\in J$. We say that $\XB$ and $\YB$ are equivalent if each one dominates the other.

Throughout this paper, we use standard notation from the theory of biorthogonal systems as can be found in \cite{HMVZ2008}. A \emph{minimal system over $\XX$} is a sequence $\XB=(x_n)_{n=1}^\infty$ in $\XX$ for which we can find another sequence $\XB^*=(x^*_n)_{n=1}^\infty$ in $\XX^*$ such that $x_n^*(x_k)=\delta_{n,k}$ for all $k$, $n\in\NN$. We call $\XB^*$ a sequence of \emph{coordinate functionals for $\XB$}. If $\XB$ is \emph{fundamental}, that is, $[\XB]=\XX$, then the sequence of coordinate functionals is unique. In this case, the \emph{support} of $x\in\XX$ relative to $\XB$ is the set
\[
\supp(x)=\enbrace{n\in\NN \tq x_n^*(x)\neq 0}.
\]
Despite the fact that a vector $x\in\spn(\XB)$ is zero if and only if $\supp(x)=\emptyset$, there are instances where the \emph{coefficient transform}
\[
\Ft\colon\XX\to\FF^\NN, \quad x\mapsto (x^*_n(x))_{n=1}^\infty
\]
is not one-to-one. If this map $\Ft$ is one-to-one, we call the minimal system $\XB$ \emph{total}. We say that a minimal system is a \emph{Markushevich} basis if it is both fundamental and total.

Any sequence $\YB=(y_n)_{n=1}^\infty$ in $\spn(\XB)\setminus\{0\}$ such that $\supp(y_k)\cap \supp(y_j)=\emptyset$ when $k\neq j$ will be called a \emph{disjoint blocking} relative to $\XB$. If we impose the stronger condition
\[
\supp(y_k)<\supp(y_{k+1})
\]
for all $k\in\NN$, then $\YB$ is a \emph{block basic sequence} relative to $\XB$. If $\YB$ is a block basic sequence, then there is a unique sequence $\alpha=(a_n)_{n=1}^\infty$ of nonzero scalars, a unique increasing sequence $\mu=(m_k)_{k=0}^\infty$ of integers with $m_0=0$, and a unique $\psi\in\Ot(\NN)$ such that 
\[
y_k=\sum_{n=1+m_{k-1}}^{m_k} a_n x_{\psi(n)}, \quad k\in\NN.
\]
If another block basic sequence $\ZB$ shares $\alpha$ and $\mu$ with $\YB$, we say that $\YB$ and $\ZB$ are \emph{equally distributed}.

Let $\XB=(x_n)_{n=1}^\infty$ be a fundamental minimal system with coordinate functionals $\XB^*=(x_n^*)_{n=1}^\infty$. Following the notation from \cite{AnsorenaMarcos2026}, given $\psi\in \It(\NN)$, we define the linear map $T_\psi=T_\psi[\XB]$ by
\[
f \mapsto \sum_{n=1}^\infty x^*_n(f)x_{\psi(n)}, \quad f\in\spn(x_n \tq n\in\NN).
\]
If $T_\psi$ extends to a bounded operator $T_\psi\colon\XX\to\XX$, we say that $\psi$ is a \emph{shifting map} relative to $\XB$, and $T_\psi$ is called a \emph{shifting operator}. The fundamental minimal system $\XB$ is \emph{spreading} if every function in $\Ot(\NN)$ is a shifting map, and the associated shifting operator is an isomorphic embedding. In turn, $\XB$ is \emph{symmetric} if every $\psi\in\Pi(\NN)$ is a \emph{shifting map}.

If $\XB$ is a spreading fundamental system of $\XX$ and $\XX$ is locally $p$-convex, then we can renorm the space so that $\XX$ is a $p$-Banach space, $\XB$ is normalized, and $T_\psi$ is an isometric embedding for all $\psi\in\Ot(\NN)$ (see \cite{AnsorenaMarcos2026}). Similarly, if $\XB$ is a symmetric fundamental system, there is an equivalent $p$-norm for $\XX$ so that $T_\psi$ is a contraction for all $\psi\in\Pi(\NN)$. Since $\Pi(\NN)$ is a group, this entails that $T_\psi$ is an isometry.

Given $A\in[\NN]_{<\infty}$, we denote by
\[
S_A=S_A[\XB]\colon\XX\to\XX, \quad f \mapsto \sum_{n\in A}x^*_n(f)\, x_n,
\]
the corresponding \emph{coordinate projection} relative to $\XB$.

A sequence $\XB=(x_n)_{n=1}^\infty$ in a quasi-Banach space $\XX$ is said to be a \emph{Schauder basis} if for every $f\in\XX$ there exists a unique sequence $(a_n)_{n=1}^\infty$ such that $f=\sum_{n=1}^\infty a_n\, x_n$. If these series converge unconditionally, we say that $\XB$ is an \emph{unconditional basis}. Any Schauder basis is a fundamental minimal system, and a fundamental minimal system $\XB$ is a Schauder basis if and only if
\[
\sup_m \norm{S_{[1\isep m]}}<\infty.
\]
If every permutation of a fundamental minimal system is a Schauder basis, then the system is an unconditional basis. This observation gives the following well-known result.

\begin{lemma}[see \cite{LinTza1977}*{Section 3}]\label{lem:ScSymUnc}
Let $\XB$ be a fundamental minimal system of a quasi-Banach space. If $\XB$ is symmetric and a Schauder basis, then it is unconditional.
\end{lemma}

For convenience, we will introduce the following notation. Given $A\in [\NN]_{<\infty}$ and a family of signs $\varepsilon=(\varepsilon_n)_{n\in A}\in\Ss_\FF^A$, we write
\[
\Ind_{\varepsilon,A}[\XB]= \sum_{n\in A} \varepsilon_n\, x_n, \quad
\Ind_A[\XB] = \sum_{n\in A} x_n .
\]
A fundamental minimal system $\XB=(x_n)_{n=1}^\infty$ in $\XX$ is \emph{democratic} if there exists a constant $C\geq 1$ such that
\[
\norm{\Ind_A[\XB]} \leq C \norm{\Ind_B[\XB]}, \quad\abs{A}\leq\abs{B}<\infty.
\]
If $\XB$ is democratic, then it is semi-normalized, that is,
\[
\sup_{n\in\NN} \norm{x_n}<\infty,\quad \inf_{n\in\NN} \norm{x_n}>0.
\]

Several properties are closely related to democracy. We say that a basis is \emph{superdemocratic} if there is a constant $C\geq 1$ such that
\[
\norm{\Ind_{\varepsilon,A}[\XB]} \leq C \norm{\Ind_{\delta,B}[\XB]}, \abs{A}\leq\abs{B}<\infty, \,\varepsilon\in \Ss_\FF^A, \, \delta\in \Ss_\FF^B.
\]
Any subsymmetric basis is superdemocratic.

We define the \emph{democracy function} of a minimal system $\XB$ in a quasi-Banach space $\XX$ as
\[
\udf[\XB](m)=\sup\enbrace{\norm{\Ind_{\varepsilon,A}[\XB]} \tq |A|\leq m, \, \varepsilon\in\Ss_\FF^A}.
\]
Note that, since $\Ind_A[\XB^*](\Ind_A[\XB])=\abs{A}$ for all $A\subset \NN$,
\[
m \le \udf[\XB](m)\udf[\XB^*](m), \quad m\in\NN.
\]
A basis is said to be \emph{bidemocratic} if a reverse inequality holds, i.e.,
\[
\udf[\XB](m)\udf[\XB^*](m)\lesssim m, \quad m\in\NN.
\]

We will write $\udf[\XB]=\udf$ when there is no ambiguity. While $\udf$ is nondecreasing by definition, an elementary convexity argument shows that if $\XX$ is a $p$-Banach space, $0<p\le 1$, then the sequence $(\udf(m)/m^{1/p})_{m=1}^\infty$ is nonincreasing (see \cite{AABW2021}*{Section 8}).

We record some well-known results for further reference.

\begin{proposition}[\cite{LinTza1977}*{Proposition 3.a.6}]\label{prop:SubSymBid}
Let $\XB$ be a subsymmetric basis of a Banach space. Then, $\XB$ is bidemocratic.
\end{proposition}

\begin{lemma}[cf.\@ \cite{AABW2021}*{Corollary 5.9}]\label{lem:FFbounded}
Let $\XB$ be an M-bounded fundamental minimal system of a quasi-Banach space. If the democracy function $\udf(m)$ is bounded and $\XB$ is bounded away from zero, then $\XB$ is equivalent to the canonical $\co$-basis.
\end{lemma}

Given a quasi-Banach space $\XX$, there is a pair $(J_\XX,\XX_c)$, where $\XX_c$ is a Banach space and $J_\XX\colon \XX \to \XX_c$ is a linear contraction such that for every Banach space $\YY$ and every linear contraction $T\colon\XX\to \YY$, there is another linear contraction $T_c\colon \XX_c \to \YY$ such that $T_c\circ J_\XX=T$. Since the pair $(J_\XX,\XX_c)$ is defined by means of a universal property, it is essentially unique. We call $\XX_c$ the \emph{Banach envelope} of $\XX$, and $J_\XX$ the \emph{envelope map}. The mapping $x^*\mapsto (x^*)_c$ is a linear isometry from $\XX^*$ to $\XX_c^*$.

The envelope map may not be one-to-one. In fact, $J_\XX$ is one-to-one if and only if $\XX$ has the \emph{point separation property}, that is, for every $x\in\XX\setminus\{0\}$ there is $x^*\in\XX^*$ such that $x^*(x)\neq 0$. The existence of a Markushevich basis implies the point separation property. If $\XB=(x_n)_{n=1}^\infty$ is a fundamental minimal system of $\XX$ with coordinate functionals $\XB^*=(x_n^*)_{n=1}^\infty$, then the envelope sequence $\XB_c=(J_\XX(x_n))_{n=1}^\infty$ is a fundamental minimal system of $\XX_c$ with coordinate functionals $((x_n^*)_c)_{n=1}^\infty$.

The \emph{bidual minimal system} $\XB^{**}$ of $\XB$ will be the sequence of coordinate functionals of $\XB^*$, regarded as a fundamental minimal system of its closed linear span. While $\XB$ dominates $\XB^{**}$, whence $\XB_c$ dominates $\XB^{**}$, there are instances where $\XB_c$ and $\XB^{**}$ are not equivalent (see \cite{AABW2021}*{Proposition 2.11}).

\begin{lemma}\label{lem:cotol1}
Let $\XB$ be a fundamental minimal system of a quasi-Banach space. Assume that $\XB$ is bounded and that $\XB^*$ is equivalent to the canonical $\co$-basis. Then, the envelope basis $\XB_c$ is equivalent to the canonical $\ell_1$-basis.
\end{lemma}

\begin{proof}
By duality, the bidual minimal system $\XB^{**}$ is equivalent to the canonical $\ell_1$-basis. Since $\XB_c$ dominates $\XB^{**}$, it dominates the canonical $\ell_1$-basis. In turn, since $\XB_c$ is a bounded sequence in a Banach space, the canonical $\ell_1$-basis dominates it.
\end{proof}

\begin{corollary}[cf.\@ \cite{DKKT2003}*{Corollary 5.7}]\label{cor:bidem-eq-l1}
Let $\XB$ be a bidemocratic fundamental minimal system of a Banach space with democracy function $\udf$. If $\udf(m)\approx m$ for $m\in\NN$, then $\XB$ is equivalent to the canonical $\ell_1$-basis.
\end{corollary}

\begin{proof}
By bidemocracy, the democracy function of the sequence $\XB^*$ of coordinate functionals is bounded. Since $\XB$ is M-bounded and dominates the bidual minimal system $\XB^{**}$, $\XB^*$ is M-bounded. Hence, combining Lemma~\ref{lem:FFbounded} with Lemma~\ref{lem:cotol1} gives the result.
\end{proof}

\begin{lemma}[\cite{AABW2021}*{Lemmas 5.6, 10.13}]\label{lem:bidem_envelope}
Let $\XB$ be a bidemocratic fundamental minimal system of a quasi-Banach space $\XX$. Then $\XB_c$ is a bidemocratic fundamental minimal system of its Banach envelope $\XX_c$. Furthermore, $\udf[\XB_c]\approx \udf[\XB]$, and $\XB^*$ is a bidemocratic minimal system of $[\XB^*]\subset\XX^*$.
\end{lemma}

The TGA, whose implementation in quasi-Banach spaces is studied in depth in \cite{AABW2021}, will also play a role in this study. Given a fundamental minimal system $\XB$ in a quasi-Banach space $\XX$ and a vector $f\in\XX$, we say that $A\in[\NN]_{<\infty}$ is a \emph{greedy set for $f$ with respect to $\XB$} if
\begin{align*}
\abs{x_k^*(f)}\leq\abs{x_n^*(f)}
\end{align*}
whenever $n\in A$ and $k\notin A$. We will denote by $\Gt[\XB](f)$ the family of all greedy sets for $f$ with respect to $\XB$. Given $m\in\NN$, the symbol $\Gt[\XB](f,m)$ stands for the family of all sets in $\Gt[\XB](f)$ whose cardinality is $m$. If the system $\XB$ is not in doubt, we write $\Gt(f)=\Gt[\XB](f)$ and $\Gt(f,m)=\Gt[\XB](f,m)$.

The behavior of bases and minimal systems with respect to greedy sets has been a fruitful field of study in modern functional analysis. Since rescaling the systems changes the greedy sets, it is customary to conduct this study on normalized, or at least semi-normalized, systems. To ensure proper implementation of the TGA, greedy sets of arbitrary cardinality must exist for all vectors. To that end, we must assume that the sequence of coordinate functionals is bounded, so that
\[
\lim_{n\to\infty} x_n^*(f)=0, \quad f\in\XX.
\]
This property ensures that $\Gt(f)$, ordered by inclusion, is a directed set. For convenience, we call the fundamental minimal systems such that
\begin{equation}
\label{eq:Wbasis}
\sup_n\max \enbrace{ \norm{x_n}, \norm{x_n^*}}<\infty
\end{equation}
\emph{W-bases}. Note that \eqref{eq:Wbasis} holds if and only if $\XB$ and $\XB^*$ are semi-normalized. Any semi-normalized Schauder basis is a W-basis.

We recall several greedy-like properties relevant to this paper. While $\XB$ is unconditional if and only if there is a constant $C$ such that
\begin{equation}\label{eq:PSB}
\norm{S_A(f)}\le C \norm{f}
\end{equation}
for all $f\in\XX$ and all $A\in[\NN]_{<\infty}$, we say that $\XB$ is \emph{quasi-greedy} if it is a W-basis for which there is a constant $C$ such that \eqref{eq:PSB} holds for all $f\in\XX$ and all $A\in\Gt(f)$. The following result was proved in its full generality in \cite{AABW2021}.
\begin{theorem}[\cite{AABW2021}*{Theorem 4.1}]\label{thm:AABWQG}
A W-basis $\XB$ of a quasi-Banach space is quasi-greedy if and only if
\[
\lim_{A\in\Gt(f)} S_A(f)=f
\]
for all $f\in\XX$.
\end{theorem}
Given $f\in\XX$, we set
\[
\varepsilon(f)=(\sgn\enpar{ x_n^*(f)})_{n=1}^\infty,
\]
where $\sgn(\lambda)=\lambda/\abs{\lambda}$ if $\lambda\neq 0$ and $\sgn(\lambda)=1$ otherwise.
We say that $\XB$ is \emph{truncation quasi-greedy} if there is a constant $C$ such that
\[
\min_{n\in A}\abs{x_n^*(f)}\norm{\Ind_{\varepsilon(f),A}}\le C \norm{f}
\]
for all $f\in\XX$ and all nonempty greedy sets $A$ of $f$. Every quasi-greedy basis is truncation quasi-greedy (see \cite{AABW2021}*{Theorem 4.13}), and the converse does not hold in general (see \cite{AABBL2023}). If the truncation quasi-greedy (or quasi-greedy, or unconditional) condition holds for vectors $f$ such that $\abs{x_n^*(f)}=1$ for all $n\in\supp(f)$, we say that $\XB$ is \emph{unconditional for constant coefficients}. A minimal system $\XB$ is superdemocratic if and only if it is democratic and unconditional for constant coefficients.

We say that $\XB$ is \emph{almost greedy} if there is a constant $C$ such that
\[
\norm{f-S_A(f)} \le C\inf\enbrace{\norm{f-S_B(f)} \tq \abs{B}\leq m}
\]
for all $m\in\NN$ and $A\in \Gt(f,m)$. If $\XB$ satisfies the more demanding condition that
\[
\norm{f-S_A(f)} \le C\inf\enbrace{\norm{f-g} \tq g\in[x_n \tq n\in B], \, \abs{B}\leq m},
\]
we say that it is \emph{greedy}.
A W-basis is almost greedy if and only if it is quasi-greedy and democratic, and it is greedy if and only if it is unconditional and democratic.

We say that the W-basis $\XB$ is \emph{semi-greedy} if there is a constant $C$ such that for every $f\in\XX$, $m\in\NN$ and $A\in\Gt(f,m)$ we can find $g\in[x_n\tq n\in A]$ such that
\[
\norm{f-g}\leq C\inf\enbrace{\norm{f-h} \tq h\in\spn(\XB)}, \,  \abs{\supp(h)}\leq m.
\]
Every almost greedy basis is semi-greedy in the general setting of quasi-Banach spaces (see \cite{BHH2025}*{Theorem 4.2}). Conversely, every semi-greedy Markushevich basis of a Banach space is almost greedy (see \cite{BL2023}*{Theorem 4.2}). Note the similarity between the latter result and the fact that symmetry and totality imply greediness by Theorem~\ref{thm:AABCOM}. Recently, Bern\'a \cite{Berna2026} found a counterexample showing that a semi-greedy Markushevich basis of a quasi-Banach space need not be almost greedy.

We will use Lorentz sequence spaces to construct suitable counterexamples witnessing Theorem~\ref{thm:BeraA} and Theorem~\ref{thm:BeraB}. As is customary, a \emph{weight} will be a sequence of positive scalars. Given a weight $\ww=(w_n)_{n=1}^\infty$, $0<p\le\infty$, and $\alpha\in\co$, we put
\begin{equation*}
\norm{\alpha}_{p,\ww}=\enpar{\sum_{m=1}^\infty \enpar{s_m \, b_m}^p \frac{w_m}{s_m}}^{1/p},
\end{equation*}
where $(b_n)_{n=1}^\infty$ is the nonincreasing rearrangement of $\abs{\alpha}$, $(s_m)_{m=1}^\infty$ is the \emph{primitive sequence} of $\ww$, defined by
\[
s_m=\sum_{n=1}^m w_n, \quad m\in\NN,
\]
and we use the standard convention if $p=\infty$. Then, we extend $\norm{\alpha}_{p,\ww}$ to $\FF^\NN$ by
\[
\norm{\alpha}_{p,\ww} =\sup\enbrace{\norm{\beta}_{p,\ww} \tq \beta\in\co, \, \abs{\beta}\le\abs{\alpha}}.
\]
The Lorentz sequence space $d_p(\ww)$ is the space consisting of all $\alpha\in\FF^\NN$ such that $\norm{\alpha}_{p,\ww}<\infty$. If the primitive sequence satisfies the doubling condition
\[
\sup_{m\in\NN} \frac{s_{2m}}{s_m}<\infty,
\]
then $d_p(\ww)$ is a quasi-Banach space (see \cite{CRS2007}). Unless $\sum_{n=1}^\infty w_n<\infty$, in which case $d_p(\ww)=\ell_\infty$, $d_p(\ww)\subset \co$. If $\sum_{n=1}^\infty w_n=\infty$ and $p<\infty$, $c_{00}$ is dense in $d_p(\ww)$. In the particular case that $\ww$ is the constant map $n\mapsto 1$, we obtain the classical Lorentz space $\ell_{1,p}$. If the weight $\ww$ is nonincreasing and $p=1$, then $\norm{\cdot}_{1,\ww}$ is a norm. Indeed, by the rearrangement inequality,
\[
\norm{\alpha}_{1,\ww}:=\sup\enbrace{\sum_{n=1}^\infty \abs{a_{\varphi(n)}} w_n \tq \varphi\in\Pi(\NN)}.
\]

\section{Preliminary results}\label{sec: preliminary}\noindent
In this section, we prove some auxiliary results that we will use later on, with a focus on totality and bidemocracy of fundamental minimal systems of quasi-Banach spaces.

\begin{lemma}\label{lemmaequivmark:new}
Let $\XB=(x_n)_{n=1}^\infty$ be a fundamental minimal system of a quasi-Banach space $\XX$ and let $(x_n^*)_{n=1}^\infty$ be its sequence of coordinate functionals. Given $x_0 \in\XX\setminus\{0\}$, the following are equivalent.
\begin{enumerate}[label=(\roman*)]
\item\label{lemmaequivmark1} $x_n^*(x_0)=0$ for all $n\in\NN$.
\item \label{lemmaequivmark2} There exists a block basic sequence $\YB=(y_n)_{n=1}^\infty$ of $\XB$ that converges to $x_0$.
\item \label{lemmaequivmark3} There exists a disjoint blocking $\YB=(y_n)_{n=1}^\infty$ of $\XB$ that converges weakly to $x_0$.
\end{enumerate}
\end{lemma}

\begin{proof}
It is clear that \ref{lemmaequivmark2} implies \ref{lemmaequivmark3}. On the other hand, if \ref{lemmaequivmark3} holds, then
\begin{align*}
x_k^*\left(x_0\right)=&\lim_{n\to \infty}x_k^*\left(y_n\right)=0
\end{align*}
for every $k\in \NN$. Thus, we obtain \ref{lemmaequivmark1}.

We now prove that \ref{lemmaequivmark1} implies \ref{lemmaequivmark2}. Let $\kappa$ be the modulus of concavity of $\XX$. Pick a null sequence $(\varepsilon_n)_{n=1}^\infty$ with $0<\varepsilon_n<\norm{x_0}$ for all $n\in\NN$. We use the fact that $x_0$ belongs to the closure of $\spn(\XB)$ to recursively construct $(y_n)_{n=1}^\infty$ in $\spn(\XB)$. Let $n\in\NN$ and assume that, if $n\ge 2$, $y_{n-1}$ has been constructed. Set $k_n=\max(\supp(y_{n-1}))$ if $n\neq 1$ and $k_n=0$ if $n=1$. Set $K_n= \norm{S_{[1 \isep k_n]}}$ and choose $z_n\in \spn(\XB)$ such that
\[
\norm{z_n-x_0}\le \frac{\varepsilon_n}{\kappa (1+K_n)}.
\]
Finally, set $y_n=z_n-S_{[1 \isep k_n]}(z_n)$ and note that $y_n\in\spn(\XB)$. For each $n\in\NN$ we have
\begin{align*}
\norm{y_n-x_0}&\le \kappa\enpar{\norm{y_n-z_n}+\norm{z_n-x_0}} \\
&= \kappa\enpar{\norm{S_{[1 \isep k_n]}(z_n-x_0)}+\norm{z_n-x_0}} \\
&\le \kappa\enpar{K_n+1} \norm{z_n-x_0}\le \varepsilon_n.
\end{align*}
Consequently, $y_n\neq 0$ for all $n\in\NN$, and $\lim_n y_n=x_0$. Since $(y_n)_{n=1}^\infty$ is a block basic sequence relative to $\XB$, we are done.

\end{proof}

\begin{lemma}\label{lemmaequivmark}
Let $\XB=(x_n)_{n=1}^\infty$ be a fundamental minimal system of a quasi-Banach space $\XX$. The following are equivalent.
\begin{enumerate}[label=(\roman*)]
\item $\XB$ is not a Markushevich basis.
\item There exists a block basic sequence $\YB=(y_n)_{n=1}^\infty$ of $\XB$ that converges to some $x_0\in\XX\setminus\enbrace{0}$.
\item There exists a disjoint blocking $\YB=(y_n)_{n=1}^\infty$ of $\XB$ that converges weakly to some $x_0\in\XX\setminus\enbrace{0}$.
\end{enumerate}
\end{lemma}

\begin{proof}
This follows immediately from Lemma~\ref{lemmaequivmark:new}.
\end{proof}

While Lemma~\ref{lemmaequivmark} is aimed at characterizing totality, the following results concern bidemocracy.

\begin{lemma}\label{lem:aux-wnullprevious}
Let $\XB=(x_n)_{n=1}^\infty$ be a fundamental minimal system of a quasi-Banach space $\XX$. If $\XB$ fails to be weakly null, then there exist $x^*\in\XX^*$, a constant $C\in(0,\infty)$ and an infinite set $A\subset\NN$ such that
\[
\sum_{n\in A} a_n \le C \norm{\sum_{n\in A} \overline{\sgn(x^*(x_n))} a_n x_n}
\]
for all eventually null nonnegative families $(a_n)_{n\in A}$.
\end{lemma}
\begin{proof}
Pick a functional $x^*$ in $\Ss_{\XX^*}$ such that $(x^*(x_n))_{n=1}^\infty$ does not converge to zero. Choose an infinite set $A\subset\NN$ with $c=\inf_{n\in A}\abs{x^*(x_n)}>0$. Given $(a_n)_{n\in A}$ as in the statement,
\begin{align*}
\norm{\sum_{n\in A} \overline{\sgn(x^*(x_n))} a_n x_n}
&\ge \abs{\sum_{n\in A} \overline{\sgn(x^*(x_n))} a_n x^*(x_n)}\\
&=\sum_{n\in A} a_n\abs{x^*(x_n)} \geq c\sum_{n\in A} a_n.\qedhere
\end{align*}
\end{proof}

\begin{lemma}\label{lem:aux-wnull}
Let $\XB=(x_n)_{n=1}^\infty$ be a fundamental minimal system of a quasi-Banach space $\XX$. If
\[
\inf_{m\in\NN}\dfrac{\udf(m)}{m}=0,
\]
then $\XB$ is weakly null.
\end{lemma}

\begin{proof}
Assume that $\XB$ is not weakly null. Let $x^*\in\XX^*$, $C\in(0,\infty)$ and $A\subset\NN$ be as in Lemma~\ref{lem:aux-wnullprevious}. Given $m\in\NN$, pick $F\subset A$ with $\abs{F}=m$. We have
\[
m \le C \norm{\sum_{n\in F} \overline{\sgn(x^*(x_n))} x_n} \le C \udf(m)
\]
and we reach a contradiction.
\end{proof}

Recall that every bounded sequence in a Banach space is dominated by the canonical $\ell_1$-basis. So, specializing Lemma~\ref{lemmabidemweaknull} below to the locally convex setting, we recover Corollary~\ref{cor:bidem-eq-l1}.

\begin{lemma}\label{lemmabidemweaknull}
Given a bidemocratic minimal system $\XB$ of a quasi-Banach space $\XX$, the following are equivalent.
\begin{enumerate}[label=(\roman*),leftmargin=*]
\item \label{lemmabidem:weaknull} $\lim_m\udf[\XB](m)/m=0$.
\item\label{lemmabidem:FailBound} $\inf_m\udf[\XB](m)/m=0$.
\item\label{lemmabidem:FailDom} $\XB$ does not dominate the canonical $\ell_1$-basis.
\item\label{lemmabidem:FailEquiv} The envelope basis $\XB_c$ is not equivalent to the canonical $\ell_1$-basis.
\item \label{lemmabidemweak1}$\XB$ is weakly null.
\item \label{lemmabidemcaucht}$\XB$ has a weakly Cauchy subsequence.
\end{enumerate}
\end{lemma}

\begin{proof}
Let $\Et$ denote the canonical $\ell_1$-basis. It is clear that \ref{lemmabidem:weaknull} implies \ref{lemmabidem:FailBound}. Since the democracy function of $\Et$ is given by $m\mapsto m$, \ref{lemmabidem:FailBound} implies \ref{lemmabidem:FailDom}. Since $\XB$ dominates $\XB_c$, \ref{lemmabidem:FailDom} implies \ref{lemmabidem:FailEquiv}. Assume \ref{lemmabidem:FailEquiv}. By Corollary~\ref{cor:bidem-eq-l1}, $\udf[\XB_c]$ is not equivalent to the mapping $m\to m$. Since $(\udf[\XB_c](m)/m)_{m=1}^\infty$ is nonincreasing,
\[
\lim_m \frac{\udf[\XB_c](m)}{m}=0.
\]
By Lemma~\ref{lem:bidem_envelope}, \ref{lemmabidem:weaknull} holds.

Assume that \ref{lemmabidem:FailEquiv} does not hold. Then, as $\Et$ is not weakly null, $\XB_c$ is not weakly null either. Since $\XX$ and $\XX_c$ have the same dual space, $\XB$ is not weakly null. So, \ref{lemmabidemweak1} implies \ref{lemmabidem:FailEquiv}.

Since \ref{lemmabidem:FailBound} implies \ref{lemmabidemweak1} by Lemma~\ref{lem:aux-wnull}, we now know that \ref{lemmabidem:weaknull}, \ref{lemmabidem:FailBound}, \ref{lemmabidem:FailDom}, \ref{lemmabidem:FailEquiv} and \ref{lemmabidemweak1} are equivalent. Assume that \ref{lemmabidemcaucht} holds, and pick $(n_j)_{j=1}^\infty$ increasing such that $(x_{n_j})_{j=1}^\infty$ is weakly Cauchy. Then, the sequence $\YB=(x_{n_{2j}}-x_{n_{2j-1}})_{j=1}^\infty$ is weakly null, and $\ZB=(x_{n_{2j}}^*)_{j=1}^\infty$ is a sequence of coordinate functionals for $\YB$. The coordinate functionals $\YB^*$ of $\YB$, regarded as a fundamental minimal system of its closed linear span, are obtained from $\ZB$ by restriction. Hence, for each $m\in\NN$,
\[
\udf[\YB](m) \udf[\YB^*](m) \le \udf[\YB](m) \udf[\ZB](m)\le \udf[\XB](2m) \udf[\XB^*](2m).
\]
So, $\YB$ inherits bidemocracy from $\XB$. Since \ref{lemmabidemweak1} implies \ref{lemmabidem:weaknull},
\[
\lim_{m\to\infty} \frac{\udf[\YB](m)}{m}=0.
\]
Since, by superdemocracy, $\udf[\XB]\lesssim \udf[\YB]$, \ref{lemmabidem:weaknull} holds.

As \ref{lemmabidemweak1} is formally stronger than \ref{lemmabidemcaucht}, the proof is complete.
\end{proof}

We close this section with a sufficient condition for bidemocracy.

\begin{lemma}\label{lem:alt_norm_bidem}
Let $\XX$ and $\YY$ be quasi-Banach spaces with fundamental minimal systems $\XB=(x_n)_{n=1}^\infty$ and $\YB=(y_n)_{n=1}^\infty$, respectively. Suppose that
\begin{enumerate}[label=(\roman*)]
\item $\YB$ is bidemocratic,
\item $\udf[\XB]\lesssim \udf[\YB]$, and
\item\label{it:3} $\XB$ dominates $\YB$.
\end{enumerate}
Then, $\XB$ is bidemocratic, and $\udf[\XB]\approx \udf[\YB]$.
\end{lemma}
\begin{proof}

Let $\XB^*$ and $\YB^*$ be the coordinate functionals of $\XB$ and $\YB$, respectively. Dualizing \ref{it:3},
we get $\udf[\XB^*]\lesssim \udf[\YB^*]$. Hence,
\[
\udf[\XB](m) \udf[\XB^*](m) \lesssim \udf[\YB](m) \udf[\YB^*](m) \lesssim m, \quad m\in \NN,
\]
and the result follows.
\end{proof}
\section{Spreading systems in quasi-Banach spaces}\label{sec:spreading-quasi-Banach}\noindent
In this section, we focus on spreading minimal systems in quasi-Banach spaces and their connections to properties such as totality, partial forms of unconditionality or greediness, and democracy.

\begin{lemma}\label{lem:weaklimit-shift-invariant}
Let $\XB=(x_n)_{n=1}^\infty$ be a spreading (resp., symmetric) fundamental minimal system of a quasi-Banach space $\XX$. If $\YB=(y_n)_{n=1}^\infty$ is a disjoint blocking of $\XB$ that converges weakly to a certain $y_0\in\XX$, then $T_\phi(y_0)=y_0$ for every $\phi\in\Ot(\NN)$ (resp., for every $\phi\in \Pi(\NN)$).
\end{lemma}

\begin{proof}
We deal with the spreading and the symmetric case simultaneously. The case $y_0=0$ is trivial, so we assume that $y_0\neq 0$. By Lemma~\ref{lemmaequivmark:new}, we may assume that $\YB$ is a block basic sequence that norm-converges to $y_0$. Using that $\phi$ is one-to-one, we recursively construct an increasing sequence $(n_k)_{k=1}^\infty$ in $\NN$ such that the sequence $(A_k)_{k=1}^\infty$ defined by
\[
A_k=\supp(y_{n_k})\cup \supp(T_\phi(y_{n_k}))
\]
satisfies $A_k<A_{k+1}$ for all $k\in\NN$. This implies that the sequence $\ZB=(z_k)_{k=1}^\infty$ defined by $z_{2k}=T_\phi(y_{n_{2k}})$ and $z_{2k-1}=y_{n_{2k-1}}$ for all $k\in\NN$ is a block basic sequence relative to $\XB$, and $(y_{n_k})_{k=1}^\infty$ and $\ZB$ are equally distributed. Let $(u_k)_{k=1}^\infty$ be a block basic sequence, and $\rho$, $\psi\in\Ot(\NN)$ be such that
\[
T_\rho(u_k)=z_k, \quad T_\psi(u_k)=y_{n_k}, \quad k\in\NN.
\]
Since symmetric fundamental systems are spreading, $T_\psi$ and $T_\rho$ are isomorphic embeddings. Since $y_{n_k}\in T_\psi(\XX)$ for all $k\in\NN$, $y_0\in  T_\psi(\XX)$. Set $u_0=T^{-1}_\psi(y_0)$. We have
\begin{align*}
T_\phi(y_0)
=&\lim_{k\to \infty}T_\phi(y_{n_{2k}})
=\lim_{k\to \infty}z_{2k},\\
y_0=&\lim_{k\to \infty}y_{n_{2k-1}}
=\lim_{k\to \infty}z_{2k-1},
\end{align*}
and, since $\enpar{T^{-1}_\psi(y_{n_k})}_{k=1}^\infty$ converges to $u_0$,
\[
T_\rho(u_0)=\lim_{k\to \infty}T_\rho\enpar{\enpar{T^{-1}_\psi(y_{n_k})}}=\lim_{k\to \infty}T_\rho(u_k)=\lim_{k\to\infty} z_k.
\]
Therefore, $T_\phi(y_0)=T_\rho(u_0)=y_0$.
\end{proof}

Theorem~\ref{thm:AABCOM} does not extend to the non-locally convex setting. Indeed, the paper \cite{AABB2026} contains an example of a symmetric Markushevich basis that is even almost greedy but is not a Schauder basis under any rearrangement. Since almost greedy bases are superdemocratic, this example rules out, in the symmetric and quasi-Banach setting, all equivalences of the Schauder basis property displayed in Theorem~\ref{thm:AABCOM} except bidemocracy. However, this does not exclude possible connections between totality, certain greedy-like properties, or weak forms of unconditionality. Our next construction answers this question in the negative for superdemocracy and for unconditionality for constant coefficients, and allows us to prove Theorem~\ref{thm:BeraB}.

\begin{proposition}\label{proposition:spread+tqg-not-markushevich}
Let $\XB=(x_n)_{n=0}^\infty$ be a spreading fundamental minimal system of a quasi-Banach space $\XX$. Assume that there exists a weight $\ww=(w_n)_{n=1}^\infty\notin\ell_1$ such that
\begin{equation}\label{eq:99}
\inf_{m\in\NN}\dfrac{\norm{\sum_{n=1}^m w_nx_n}}{\sum_{n=1}^m w_n} = 0.
\end{equation}
If for each $n\in\NN$ we define
\[
y_n = x_n+x_0,
\]
and set $\YB=(y_n)_{n=1}^\infty$ and $\YB^*=(x^*_n)_{n=1}^\infty$, then the following hold.
\begin{enumerate}[label=(\roman*)]
\item $\YB$ is a fundamental minimal system of $\XX$ with biorthogonal functionals $\YB^*$.
\item\label{lem:spread+tqg-not-markushevich:spread} $\YB$ is spreading.
\item\label{lem:spread+tqg-not-markushevich:sym} If $\XB$ is symmetric, so is $\YB$.
\item $\YB$ fails to be a Markushevich basis.
\item\label{lem:spread+tqg-not-markushevich:semigreedy} If $\XB$ is unconditional, then $\YB$ is semi-greedy.
\item \label{lem: spread:tqg} If $\XB$ is truncation quasi-greedy and
\[
K:=\inf_{m\in\NN} \dfrac{\udf[\XB](m)}{m}>0,
\]
then $\YB$ is truncation quasi-greedy and superdemocratic, with $\udf[\XB]\approx\udf[\YB]$.
\end{enumerate}
\end{proposition}
\begin{proof}
The sequence $\YB$ inherits symmetry or spreadingness from $\XB$, so \ref{lem:spread+tqg-not-markushevich:spread} and \ref{lem:spread+tqg-not-markushevich:sym} hold. Since
\[
\frac{\sum_{n=1}^mw_ny_n}{\sum_{n=1}^mw_n}-x_0=\frac{\sum_{n=1}^mw_n x_n}{\sum_{n=1}^mw_n},
\]
$\YB$ is fundamental by \eqref{eq:99}. Besides, $y_n^*(x_0)=0$ for all $n\in\NN$, so $\YB$ is not total.

To prove \ref{lem:spread+tqg-not-markushevich:semigreedy}, we follow some ideas from \cite{Berasategui2024}*{Example 1.2.16}. By the Aoki--Rolewicz theorem, we may assume without loss of generality that $\XX$ is a $p$-Banach space for some $0<p\le 1$. Additionally, using a suitable renorming, we may assume that $\XB$ is normalized and isometrically subsymmetric, that the shifting operator $T_\psi[\XB]\colon\XX\to\XX$ is an isometric embedding for every $\psi\in\Ot(\ZZ_{\ge 0})$, and that the coordinate projection $S_A[\XB]\colon\XX\to\XX$ is a contraction for any $A\in[\ZZ_{\ge 0}]_{<\infty}$. The latter property implies that its complementary projection 
\[
S_A^c[\XB]:=\Id_\XX-S_A[\XB]
\]
is also a contraction.

Pick $f\in\XX$, $m\in\NN$, $B\in\Gt[\YB](f,m)$, and $A\subset\NN$ with $\abs{A}=m$. We must prove that for every $h\in [y_n \tq n\in A]$ there is $g\in [y_n \tq n\in B]$ such that $\norm{f-g}\leq C\norm{f-h}$, where the constant $C\in(0,\infty)$ only depends on $\YB$. Let $\pi\colon B\setminus A\to A\setminus B$ be an order-preserving bijection. Set $h=\sum_{n\in A} a_n\, y_n$. Define $(b_n)_{n\in B}$ as
\[
b_n =
\begin{cases}
a_n, \quad\text{if $n\in A$},\\
a_{\pi(n)}\quad\text{if $n\in B\setminus A$},
\end{cases}
\]
and set $g=\sum_{n\in B}  b_n \, y_n$. We compare the coefficients of $f-g$ and $f-h$ relative to $\XB$ by distinguishing between several cases.

\begin{enumerate}[label=(\alph*),leftmargin=*,widest=a]
\item\label{case:A} If $k\in\NN \setminus (A\cup B)$, then $x_k^*(f-g)=x_k^*(f-h)=y_k^*(f)$.
\item\label{case:B} If $k\in A\cap B$, then $x_k^*(f-g)=x_k^*(f-h)=x_k^*(f)-a_k$.
\item\label{case:C} If $k=0$, then
\[
x_k^*(f-g)=x_0^*(f)-\sum_{n\in A} a_n=x_0^*(f)-\sum_{n\in B} b_n=x_k^*(f-h).
\]
\item\label{case:D} If $k\in B\setminus A$, then $x_k^*(f-g)=x_k^*(f)-a_{\pi(k)}$ and $x_k^*(f-h)=x_k^*(f)$.
\item\label{case:E} If $k\in A\setminus B$, then $x_k^*(f-g)=x_k^*(f)$ and $x_k^*(f-h)=x_k^*(f)-a_k$.
\end{enumerate}

Pick $k\in B\setminus A$ and set $\mu_k=\abs{x_k^*(f)}+\abs{a_{\pi(k)}}$. By \ref{case:D} and \ref{case:E},
\begin{align*}
&\max\enbrace{ \abs{x_k^*(f-g)}, \abs{x_{\pi(k)}^*(f-g)}}
\le \max\enbrace{ \abs{x_k^*(f)}+\abs{a_{\pi(k)}} , \abs{x_{\pi(k)}^*(f) }} \\
&\le \max\enbrace{ \abs{x_k^*(f)}+\abs{a_{\pi(k)}} , \abs{x_k^*(f) }}
= \mu_k
\end{align*}
and, since $\pi(k)$ is outside $B$ while $k$ belongs to $B$,
\begin{align*}
\abs{x_k^*(f-h)}+\abs{x_{\pi(k)}^*(f-h)}
&\ge\abs{x_k^*(f-h)}+\frac{1}{2}\abs{x_{\pi(k)}^*(f-h)}\\
&\ge\abs{x_k^*(f)}+\frac{1}{2}\enpar{\abs{a_{\pi(k)}} - \abs{x_{\pi(k)}^*(f)}}\\
&\ge\abs{x_k^*(f)}+\frac{1}{2}\enpar{\abs{a_{\pi(k)}} - \abs{x_k^*(f)}}
=\frac{1}{2} \mu_k.
\end{align*}
Set $\mu=\sum_{k\in B\setminus A} \mu_k \, x_k$. We have
\begin{align*}
\norm{\mu}^p & \le 2^p \norm{ \sum_{k\in B\setminus A} \enpar{\abs{x_k^*(f-h)}+\abs{x_{\pi(k)}^*(f-h)} } x_k }^p\\
&\le 2^p \enpar{ \norm{ \sum_{k\in B\setminus A} \abs{x_k^*(f-h)} x_k }^p +\norm{ \sum_{k\in B\setminus A} \abs{x_{\pi(k)}^*(f-h)} x_k }^p }\\
&= 2^p \enpar{ \norm{ \sum_{k\in B\setminus A} \abs{x_k^*(f-h)} x_k}^p +\norm{ \sum_{k\in B\setminus A} \abs{x_{\pi(k)}^*(f-h)} x_{\pi(k)} }^p }\\
&\le 2^{p+1}\norm{f-h}^p,
\end{align*}
and, if $D=A\triangle B$,
\begin{align*}
\norm{S_D[\XB](f-g)}^p
&\le \norm{ \sum_{k\in B\setminus A} \mu_k\, x_k + \sum_{k\in B\setminus A} \mu_k\, x_{\pi(k)}}^p\\
&\le \norm{ \sum_{k\in B\setminus A} \mu_k\, x_k}^p + \norm{\sum_{k\in B\setminus A} \mu_k\, x_{\pi(k)}}^p=2\norm{\mu}^p.
\end{align*}
By \ref{case:A}, \ref{case:B}, and \ref{case:C},
\[
\norm{S_D^c[\XB](f-g)}=\norm{S_D^c[\XB](f-h)}\le \norm{f-h}.
\]
Summing up,
\[
\norm{f-g}^p \le \norm{ S_D[\XB](f-g)}^p +\norm{ S_D^c[\XB](f-g)}^p \le \enpar{2^{p+2}+1} \norm{f-h}^p.
\]

Finally, we prove \ref{lem: spread:tqg}. Notice that combining truncation quasi-greediness and spreadingness yields superdemocracy. Let $\kappa$ be the modulus of concavity of $\XX$, and let $K_t$ and $K_s$ be the truncation quasi-greedy and superdemocracy constants of $\XB$, respectively. Also, note that there is a constant $C_1$ such that
\[
\norm{f-x_0^*(f) x_0} \le C_1 \norm{f}, \quad f\in\XX.
\]
Indeed, we could use $C_1=\kappa(1+\norm{x_0^*} \norm{x_0})$. Fix $A\subset\NN$ and $\varepsilon=(\varepsilon_n)_{n\in A}\in \Ss_\FF^A$. Given $f\in\XX$ such that $\abs{x_n^*(f)}\ge 1$ for all $n\in A$ and $\abs{x_n^*(f)}\le 1$ for all $n\in \NN\setminus A$,
\begin{equation*}
\norm{\Ind_{\varepsilon,A}[\XB]} \le K_s K_t \norm{f-x_0^*(f) x_0}\le C_1 K_s K_t \norm{f}.
\end{equation*}
If $g=\Ind_{\varepsilon,A}[\YB]$ and $s=\sum_{n\in A} \varepsilon_n$, then $g=\Ind_{\varepsilon,A}[\XB] + s x_0$. Hence,
\[
\norm{g}\leq\kappa\enpar{\norm{\Ind_{\varepsilon,A}[\XB]}+\abs{s} \norm{x_0}}
\leq\kappa\enpar{1+\frac{K_s}{K}\norm{x_0}}\norm{\Ind_{\varepsilon,A}[\XB]}.
\]
Conversely,
\[
\norm{\Ind_{\varepsilon,A}[\XB]}=\norm{g- x_0^*(g) x_0}\le C_1 \norm{g}.
\]
Hence, $\YB$ is truncation quasi-greedy and inherits superdemocracy from $\XB$. Furthermore, $\udf[\XB]\approx\udf[\YB]$.
\end{proof}

Now we are ready to prove the main result of this section.

\begin{proof}[Completion of the proof of Theorem~\ref{thm:BeraB}]
Choose $1<p\le \infty$, and let $\Et=(e_n)_{n=1}^\infty$ be the canonical $\ell_{1,p}$-basis. We have
\[
\norm{\sum_{n=1}^m e_{n+1}}=\enpar{\sum_{n=1}^m n^{p-1}}^{1/p} \approx m, \quad m\in\NN,
\]
and
\[
\norm{\sum_{n=1}^m \frac{1}{n} e_{n+1}}=\enpar{\sum_{n=1}^m \frac{1}{n}}^{1/p}, \quad m\in\NN,
\]
with the obvious modification if $p=\infty$. Therefore, choosing $\ww=(1/n)_{n=1}^\infty$, Proposition~\ref{proposition:spread+tqg-not-markushevich} applies to $\XB=(e_{n+1})_{n=0}^\infty$.
\end{proof}
\section{Spreading systems in Banach spaces}\label{sec:spreading-Banach}\noindent
In this section, in addition to proving Theorem~\ref{thm:dem-spread-equivalences}, we discuss whether the properties of minimal systems involved in its statement are equivalent for spreading fundamental minimal systems in Banach spaces. We also characterize all minimal systems that are spreading, but fail to be Markushevich bases.

We start by recording a well-known consequence of the Hahn--Banach theorem for further reference.

\begin{lemma}\label{lem:weaknull-span}
Let $\XB=(x_j)_{j\in J}$ be a net in a Banach space $\XX$, and $x\in\XX$. Let $\YY$ be a closed subspace of $\XX$ such that $x_j\in\YY$ for all $j\in J$. If $\XB$ converges weakly to $x$, then $x\in\YY$ and $\XB$ converges to $x$ in the weak topology of $\YY$.
\end{lemma}

Lemma~\ref{lem:weaknull-span} yields the following result. We emphasize that, since the proof uses the Hahn--Banach theorem, this argument does not apply in non-locally convex spaces.

\begin{lemma}\label{lem:weak-spread-convergence}
Let $\XB=(x_n)_{n=1}^\infty$ be a spreading minimal system of a Banach space $\XX$, and let $\Ut=(u_n)_{n=1}^\infty$
and $\ZB=(z_n)_{n=1}^\infty$ be subsequences of $\XB$. Then, $\Ut$ converges weakly to the origin if and only if $0\in \overline{\enbrace{z_n\tq n\in\NN}}^w$.
\end{lemma}

\begin{proof}
Assume that $\Ut$ does not converge weakly to the origin. Then, there is a subsequence $\YB=(y_n)_{n=1}^\infty$ of $\Ut$ such that
\[
0\notin \overline{\enbrace{y_n\tq n\in\NN}}^w.
\]
By Lemma~\ref{lem:weaknull-span}, the origin does not belong to the weak closure of $\YB$ in $[\YB]$. Since $\YB$ and $\ZB$ are equivalent, the origin does not belong to the weak closure of $\ZB$ in $[\ZB]$. By Lemma~\ref{lem:weaknull-span}, the origin is not a weak cluster point of $\ZB$.

Conversely, if $\lim_n u_n=0$ weakly, then, by Lemma~\ref{lem:weaknull-span}, $\Ut$ converges to the origin in the weak topology of $[u_n \tq n \in\NN]$. Since $\Ut$ and $\ZB$ are equivalent, $\ZB$ converges weakly to the origin in $[\ZB]$. This implies that the origin is a cluster point of $\ZB$ in the weak topology of $[\ZB]$. Hence, the origin is a weak cluster point of $\ZB$ in $\XX$.
\end{proof}

In the following results, we use standard facts from the theory of spreading models. The reader can find background on this theory in \cite{AnsorenaMarcos2026}.

\begin{lemma}\label{lem:spread+wnull=>unc}
Let $\XB=(x_n)_{n=1}^\infty$ be a spreading fundamental minimal system of a Banach space $\XX$. If the origin is in the weak closure of $\XB$, then $\XB$ is unconditional.
\end{lemma}
\begin{proof}
Since $\XB$ is spreading, it is semi-normalized. There exists a suitable subsequence $\YB$ of $\XB$ that has a spreading model (see \cite{BeauzamyLapreste1983}*{Chapitre I, Proposition 1}). By Lemma~\ref{lem:weak-spread-convergence}, $\YB$ converges weakly to zero. Hence, the unit vector system $\St$ of the spreading model of $\YB$ is unconditional (see \cite{BeauzamyLapreste1983}*{Chapitre I, Proposition 4}). Besides, by spreadingness, $\XB$ is equivalent to $\St$, and unconditionality follows.
\end{proof}

\begin{lemma}\label{lem:spread+ucc+weak-cauchy=>subsym}
Let $\XX$ be a Banach space and $\XB=(x_n)_{n=1}^\infty$ be a spreading fundamental minimal system in $\XX$. Assume that $\XB$ is unconditional for constant coefficients and has a weakly Cauchy subsequence. Then, $\XB$ is weakly null and unconditional.
\end{lemma}
\begin{proof}
We denote by $\ZB = (z_k)_{k=1}^\infty$ the weakly Cauchy subsequence of $\XB$ and set $\YB=(z_{2k-1}-z_{2k})_{k=1}^\infty$. This new sequence is a minimal system that is weakly null and inherits spreadingness from $\XB$. Hence, by Lemma~\ref{lem:weaknull-span} and Lemma~\ref{lem:spread+wnull=>unc}, $\YB$ is unconditional. By Proposition~\ref{prop:SubSymBid}, $\YB$ is also bidemocratic, so, by Lemma~\ref{lemmabidemweaknull}, 
\[
\lim_{m\to\infty} \frac{\udf[\YB](m)}{m}=0.
\]
Since $\ZB$ is unconditional for constant coefficients,
\[
\udf[\YB](m) \approx \udf[\ZB](2m), \quad  m\in\NN.
\]
As $\XB$ and $\ZB$ are equivalent, this gives 
\[
\lim_{m\to\infty} \frac{\udf[\XB](m)}{m}=0.
\]
It follows from Lemma~\ref{lem:aux-wnull} that $\XB$ is weakly null and, from Lemma~\ref{lem:spread+wnull=>unc}, that $\XB$ is unconditional.
\end{proof}

We are now ready to address the proof of the main result of this section, followed by that of Proposition~\ref{prop:spread_not_markushevich}.

\begin{proof}[Completion of the proof of Theorem~\ref{thm:dem-spread-equivalences}]
Condition \ref{cor:dem-spread-equivalences:unc} implies \ref{cor:dem-spread-equivalences:bidem} by a direct application of Proposition~\ref{prop:SubSymBid}. It is clear that \ref{cor:dem-spread-equivalences:bidem} implies \ref{cor:dem-spread-equivalences:superdem}, and that \ref{cor:dem-spread-equivalences:superdem} implies \ref{cor:dem-spread-equivalences:ucc}. So, it suffices to prove that \ref{cor:dem-spread-equivalences:ucc} implies \ref{cor:dem-spread-equivalences:unc}. To that end, we recall that, by Rosenthal's $\ell_1$-theorem \cite{Rosenthal1974}, either $\XB$ has a weakly Cauchy subsequence, or, by spreadingness, it is equivalent to the canonical $\ell_1$-basis. Applying Lemma~\ref{lem:spread+ucc+weak-cauchy=>subsym} yields the desired connection.
\end{proof}

\begin{proof}[Completion of the proof of Proposition~\ref{prop:spread_not_markushevich}]
Since $\XB$ is spreading, it is bounded (see \cite{AnsorenaMarcos2026}*{Lemma 2.4}) and has no Schauder subsequences. Hence, by \cite{AlbiacKalton2016}*{Theorem 1.5.6}, $\XB$ is relatively weakly compact, and $0$ is not a weak accumulation point of $\XB$. By the Eberlein--\v{S}mulian theorem (see, e.g., \cite{AlbiacKalton2016}*{Theorem 1.6.3}), there is a subsequence $\Ut=(u_n)_{n=1}^\infty$ of $\XB$ that converges weakly to a certain nonzero element $u_0$. By Lemma~\ref{lem:weaknull-span}, $u_0\in [u_n \tq n\in \NN]$ and $\Ut$ converges to $u_0$ in the weak topology of $[u_n \tq n\in \NN]$. Since $\Ut$ and $\XB$ are equivalent, $\XB$ converges weakly to a certain nonzero element $y_0$.

Now we move on to the sequence $\YB=(y_n)_{n=1}^\infty$ given by $y_n=x_n-y_0$ for all $n\in\NN$. If $(x_n^*)_{n=1}^\infty$ is the sequence of coordinate functionals of $\XB$, then, for all $k\in\NN$,
\[
x_k^*(y_0)=\lim_{n\to\infty} x_k^*(x_n)=\lim_{n\to\infty} \delta_{k,n}=0.
\]
Hence, $\YB$ is a fundamental minimal system of its closed linear span $\YY=[\YB]$ with coordinate functionals $(y_n^*)_{n=1}^\infty$ given by $y_n^*=\restr{x_n^*}{\YY}$ for all $n\in\NN$. Since $\YB$ is weakly null, by Lemma~\ref{lem:weaknull-span}, $\YB$ converges to zero in the weak topology of $\YY$.

By assumption, for each $\phi\in\Ot(\NN)$ there is an isomorphic embedding $T_\phi\colon\XX\to\XX$ such that $T_\phi(x_n)=x_{\phi(n)}$ for all $n\in\NN$. By Lemma~\ref{lem:weaklimit-shift-invariant}, $T_\phi(y_n)=y_{\phi(n)}$ for all $n\in\NN$. Hence, $\phi$ is a shifting map relative to $\YB$, and the associated shifting operator is $\restr{T_\phi}{\YY}$, which is an isomorphic embedding. We have proved that $\YB$ is spreading. It follows from Lemma~\ref{lem:spread+wnull=>unc} that $\YB$ is an unconditional basis of $\YY$.

To prove that by adding the vector $y_0$ to $\YB$ we obtain a subsymmetric basis of $\XX$, it suffices to prove that $y_0\notin\YY$. Assume, towards a contradiction, that this is not the case. Then, there is a series expansion
\[
y_0=\sum_{n=1}^\infty a_n \, y_n.
\]
Given $k\in\NN$, let $\psi\in\Ot(\NN)$ be defined by $\psi(\NN)=\NN\setminus\{k\}$. By Lemma~\ref{lem:weaklimit-shift-invariant},
\[
\sum_{n=1}^\infty a_n \, y_n=y_0=T_\psi(y_0)=\sum_{n=1}^\infty a_n \, y_{\psi(n)}.
\]
By comparing coefficients, we deduce that $a_k= 0$. Consequently, $y_0=0$. This ends the proof.
\end{proof}

We close this section with some results linking spreading and semi-greedy minimal systems. To contextualize the first result, it is worth noting that, by Lemma~\ref{lem:aux-wnullprevious}, a subsymmetric basis of a Banach space fails to be weakly null if and only if it is equivalent to the canonical $\ell_1$-basis.

\begin{corollary}\label{prop:char_spread_not_markushevich}
Let $\XB=(x_n)_{n=0}^\infty$ be a subsymmetric (resp., unconditional and symmetric) basis of a Banach space $\XX$ that is not equivalent to the canonical $\ell_1$-basis. Let $(x_n^*)_{n=0}^\infty$ be the coordinate functionals of $\XB$. Define for each $n\in\NN$
\[
y_n = x_n+x_0.
\]
Then, $\YB=(y_n)_{n=1}^\infty$ is a fundamental minimal system of $\XX$ with coordinate functionals $\YB^*=(x_n^*)_{n=1}^\infty$. Furthermore, $\YB$ satisfies the following properties.
\begin{enumerate}[label=(\roman*)]
\item $\YB$ is spreading (resp., symmetric).
\item $\YB$ is semi-greedy.
\item $\YB$ fails to be a Markushevich basis.
\end{enumerate}
\end{corollary}

\begin{proof}
By Proposition~\ref{prop:SubSymBid} and Corollary~\ref{cor:bidem-eq-l1}, $\udf[\XB]$ is not equivalent to the map $m\mapsto m$. Since $(\udf[\XB](m)/m)_{m=1}^\infty$ is nonincreasing,
\[
\lim_m \frac{\udf[\XB](m)}{m}=0.
\]
Hence, Proposition~\ref{proposition:spread+tqg-not-markushevich} applies with $(w_n)_{n=1}^\infty$ given by $w_n=1$ for all $n\in\NN$.
\end{proof}

\begin{theorem}
Let $\XB=(x_n)_{n=1}^\infty$ be a fundamental minimal system of a Banach space. If $\XB$ is spreading and fails to be total, then it is semi-greedy.
\end{theorem}

\begin{proof}
By Proposition~\ref{prop:spread_not_markushevich}, there is $y_0\in\XX$ such that the sequence $\YB=(y_n)_{n=0}^\infty$ given by $y_n=x_n-y_0$ for all $n\in\NN$ is a subsymmetric basis of $\XX$. Furthermore, $\YB$ is weakly null, whence it is not equivalent to the canonical $\ell_1$-basis. By Corollary~\ref{prop:char_spread_not_markushevich}, $\XB$ is semi-greedy.
\end{proof}
\section{Symmetric systems in quasi-Banach spaces}\label{sec:sym-qBanach}\noindent

In this section, we study bidemocracy and its interplay with symmetry and other properties of interest. We will show that symmetric bidemocratic minimal systems need not be unconditional. In fact, there are quasi-Banach spaces with Markushevich bases that are symmetric, bidemocratic, and almost greedy but fail to be Schauder bases under any reordering.

For the first of these results, it will be convenient to extend several notions introduced in Section~\ref{sec:Background} for fundamental minimal systems of quasi-Banach spaces, which also make sense for the unit vector system $\Et=(\ee_n)_{n=1}^\infty$ of the vector space $\FF^\NN$, regardless of whether we are considering a quasi-Banach space $\XX$ with $\coo\subset\XX\subset\FF^\NN$. Specifically, given $\psi\in \It(\NN)$, we denote by $T_\psi\colon\FF^\NN \to \FF^\NN$ the shifting operator relative to $\Et$, that is, if $\alpha=(a_n)_{n=1}^\infty$, then $T_\psi(\alpha)=(b_n)_{n=1}^\infty$, where $b_k=a_n$ if $k=\psi(n)$ for some $n\in\NN$, and $b_k=0$ otherwise. Similarly, $S_A\colon\FF^\NN\to \FF^\NN$ denotes the coordinate projection associated with $A\subset\NN$. Given a sequence $\alpha=(a_n)_{n=1}^\infty$ of scalars, we set $\varepsilon(\alpha)= (\sgn (a_n))_{n=1}^\infty$ and $\supp(\alpha)=\enbrace{n\in\NN \tq a_n\neq 0}$. We say that $A\in [\NN]_{<\infty}$ is a greedy set for $\alpha$ if $\abs{a_n}\ge \abs{a_k}$ for all $n\in A$ and $k\in\NN\setminus A$.

Given $A\subset\NN$ and $\varepsilon\in \Ss_\FF^A$, we denote by $\Ind_A$ and $\Ind_{\varepsilon,A}$ the corresponding indicator function and signed indicator function relative to $\Et$. Similarly, if $\abs{A}<\infty$, $\Ind_A^*$ stands for the indicator function relative to the coordinate functionals $(\ee_n^*)_{n=1}^\infty$ of the unit vector system of $\FF^\NN$, that is,
\[
\Ind_A^*\colon \FF^\NN \to \FF, \quad (a_n)_{n=1}^\infty \mapsto \sum_{n\in A} a_n.
\]

Given a nonincreasing weight $\ww=(w_n)_{n=1}^\infty$ with primitive sequence $(s_m)_{m=1}^\infty$, the \emph{Marcinkiewicz sequence space} $m(\ww)$ is the Banach space consisting of all scalar-valued sequences $\alpha=(a_n)_{n=1}^\infty$ such that
\[
\norm{\alpha}_{\cM, \ww}:=\sup_{\substack{m\in\NN\\ \abs{A}=m}} \frac{\sum_{n\in A} \abs{a_n}}{s_m} <\infty.
\]
Unless $\lim_n w_n>0$, in which case $m(\ww)=\ell_\infty$, $m(\ww)\subset \co$. The \emph{little Marcinkiewicz sequence space} $m_0(\ww)$ is the closure of $\coo$ in $m(\ww)$. Unless $\sum_{n=1}^\infty w_n<\infty$, in which case $m_0(\ww)=m(\ww)=\ell_1$, the sequence $\alpha$ lies in $m_0(\ww)$ if and only if
\begin{equation*}
\lim_{m\to\infty}\sup_{\abs{A}=m}\frac{\sum_{n\in A} \abs{a_n}}{s_m}=0.
\end{equation*}
If $\alpha\in\co$ and $(b_n)_{n=1}^\infty$ is the nonincreasing rearrangement of $\abs{\alpha}$,
\[
\norm{\alpha}_{\cM,\ww}=\sup_{m\in\NN} \frac{\sum_{n=1}^m b_n}{s_m}.
\]

For every $\alpha=(a_n)_{n=1}^\infty\in\co$ we define
\[
G_m(\alpha)= \sup_{A\in\Gt(\alpha,m)} \abs{\Ind_A^*(\alpha)}, \quad m\in\NN,
\]
and
\[
\norm{\alpha}_{\cG, \ww}=\sum_{m=1}^\infty \frac{w_m}{m} G_m(\alpha).
\]
The functional $\norm{\cdot}_{\cG,\ww}$ is a Lorentz norm on nonnegative sequences. Indeed, if $\alpha\in\co$ and $(b_n)_{n=1}^\infty$ is the nonincreasing rearrangement of $\abs{\alpha}$, by the Abel summation formula,
\begin{equation}\label{eq:Gpos}
\norm{\alpha}_{\cG, \ww}\le \norm{\abs{\alpha}}_{\cG, \ww} =\sum_{m=1}^\infty \frac{w_m}{m} \sum_{n=1}^m b_n=\sum_{n=1}^\infty \enpar{\sum_{m=n}^\infty \frac{w_m}{m}} b_n.
\end{equation}

The functional $\norm{\cdot}_{\cG,\ww}$, as well as $\norm{\cdot}_{\cM,\ww}$ and $\norm{\cdot}_{1,\ww}$, is positive definite and absolutely homogeneous. The following result is a substitute for the triangle inequality.

\begin{lemma}\label{lem:Anso1}
Let $\ww$ be a nonincreasing weight. Then, for all $\alpha$, $\beta\in\co$,
\begin{equation*}
\norm{\alpha+\beta}_{\cG, \ww}
\le \norm{\alpha}_{\cG, \ww}+\norm{\beta}_{\cG, \ww}+2\enpar{\norm{\alpha}_{1,\ww}+\norm{\beta}_{1,\ww}+\norm{\alpha+\beta}_{1,\ww}}.
\end{equation*}
\end{lemma}
\begin{proof}
Pick $m\in\NN$. Let $(a_n)_{n=1}^\infty$, $(b_n)_{n=1}^\infty$, and $(d_n)_{n=1}^\infty$ be the nonincreasing rearrangements of $\abs{\alpha}$, $\abs{\beta}$ and $\abs{\alpha+\beta}$, respectively. Pick $A\in\Gt(\alpha,m)$ and $B\in\Gt(\beta,m)$ and $D\in \Gt(\alpha+\beta,m)$.
The expansion
\begin{multline*}
\Ind_D^*(\alpha+\beta)-\Ind_A^*(\alpha)-\Ind_B^*(\beta)\\
=\Ind^*_{(B\cup D)\setminus A}(\alpha)+\Ind^*_{(A\cup D)\setminus B}(\beta)-\Ind^*_{(A\cup B)\setminus D}(\alpha+\beta)
\end{multline*}
yields
\[
G_m(\alpha+\beta) \le G_m(\alpha) +G_m(\beta) + 2 m\enpar{a_{m+1}+b_{m+1}+d_{m+1}}.
\]
Summing over $m$ with weights $w_m/m$, we obtain
\[
\norm{\alpha+\beta}_{\cG, \ww}\le \norm{\alpha}_{\cG, \ww}+\norm{\beta}_{\cG, \ww}+2\sum_{m=1}^\infty w_m \enpar{a_{m+1}+b_{m+1}+d_{m+1}}.
\]
This estimate readily implies the desired inequality.
\end{proof}

\begin{lemma}\label{lem:Anso2}
Let $\ww$ be a nonincreasing weight. Let $(\alpha_j)_{j=1}^\infty$ be a sequence in $\co$, and let $\alpha\in\co$ be such that $\lim_j \alpha_j=\alpha$ uniformly. Then,
\begin{equation*}
\norm{\alpha}_{\cG,\ww} \le\liminf_j \norm{\alpha_j}_{\cG,\ww}+ 2 \liminf_j \norm{\alpha_j}_{1,\ww}.
\end{equation*}
\end{lemma}

\begin{proof}
Let $\beta$ and $\beta_j$, $j\in\NN$, be the nonincreasing rearrangements of $\abs{\alpha}$ and $\abs{\alpha_j}$, respectively. Since
\[
\beta\le \liminf_j \beta_j
\]
pointwise, by Fatou's lemma,
\begin{equation}\label{eq:AnsoLorentz}
\norm{\alpha}_{1,\ww}\le \liminf_j \norm{\alpha_j}_{1,\ww}.
\end{equation}
Set $\beta=(b_k)_{k=1}^\infty$. For each $k\in \NN$, choose $A_k\in \Gt(\alpha,k)$ so that
\[
\abs{\Ind_{A_k}^*(\alpha)}=G_k(\alpha).
\]
Also, for each $j,k\in \NN$, choose $A_{k,j}\in \Gt(\alpha_j,k)$ so that
\[
\abs{\Ind_{A_{k,j}}^*(\alpha_j)}=G_k(\alpha_j).
\]
For every $k\in \NN$ we have

\begin{align*}
G_k(\alpha)=&\abs{\Ind_{A_k}^*(\alpha-\alpha_j)+\Ind_{A_{k,j}}^*(\alpha_j)+\Ind^*_{A_k\setminus A_{k,j}}(\alpha_j)-\Ind_{A_{k,j}\setminus A_k}^*(\alpha_j)}\\
\le& \abs{\Ind_{A_k}^*(\alpha-\alpha_j)}+\abs{\Ind_{A_{k,j}}^*(\alpha_j)}+\abs{\Ind^*_{A_k\setminus A_{k,j}}(\alpha_j)-\Ind_{A_{k,j}\setminus A_k}^*(\alpha_j)}\\
\le& \abs{\Ind_{A_k}^*(\alpha-\alpha_j)}+G_k(\alpha_j)+ \abs{A_k\triangle A_{k,j}} \norm{S_{A_k\triangle A_{k,j}}(\alpha_j)}_\infty\\
=& \abs{\Ind_{A_k}^*(\alpha-\alpha_j)}+G_k(\alpha_j)+ \abs{A_k\triangle A_{k,j}} \norm{S_{A_{k,j}\setminus A_k}(\alpha_j)}_\infty\\
\le& \abs{\Ind_{A_k}^*(\alpha-\alpha_j)}+G_k(\alpha_j)+2k \norm{S_{\NN \setminus A_k}(\alpha_j)}_\infty.
\end{align*}

By uniform convergence, taking the limit inferior on both sides yields
\[
G_k(\alpha)\le \liminf_j G_k(\alpha_j)+ 2k \norm{S_{\NN \setminus A_k}(\alpha)}_\infty \le \liminf_j G_k(\alpha_j)+ 2k b_k.
\]
Summing over $k$ with weights $w_k/k$ and using Fatou's lemma, we get
\[
\norm{\alpha}_{\cG,\ww}\le 2\norm{\alpha}_{1,\ww}+\liminf_j \norm{\alpha_j}_{\cG, \ww}.
\]
Combining this estimate with \eqref{eq:AnsoLorentz} ends the proof.
\end{proof}

Given two nonincreasing weights $\uu$ and $\ww$, we consider the space
\[
\BB_{\uu,\ww}=\enbrace{\alpha\in \co \tq \norm{\alpha}_{\uu,\ww}:=\max\enbrace{ \norm{\alpha}_{\cM, \uu}, \norm{\alpha}_{1,\ww}, \norm{\alpha}_{\cG,\ww}}<\infty}.
\]

 When $\BB_{\uu,\ww}$ is a quasi-Banach space containing the unit vector system, we denote by $\BB_{\uu,\ww}^0$ the closure of $\coo$ in $\BB_{\uu,\ww}$. We give conditions ensuring that this situation occurs. Some properties of the space $\BB_{\uu,\ww}$ depend on the pair $(\uu,\ww)$ satisfying a two-weight Dini condition.

\begin{definition}
Let $\uu$ and $\ww=(w_n)_{n=1}^\infty$ be nonincreasing weights, and let $(s_m)_{m=1}^\infty$ be the primitive sequence of $\uu$. Given $0<K<\infty$, we say that $(\uu,\ww)$ is a head Dini pair with constant $K$ if
\[
\frac{s_m}{m}\sum_{n=1}^m w_n \le K, \quad m\in \NN.
\]
We say that $(\uu,\ww)$ is a tail Dini pair with constant $K$ if
\[
s_m\sum_{n=m+1}^\infty\frac{w_n}{n}\le K, \quad m\in \NN.
\]
If $(\uu,\ww)$ is both a head Dini pair and a tail Dini pair with constant $K$, we say that $(\uu,\ww)$ is a \emph{two-sided Dini pair} with constant $K$.
\end{definition}

This terminology is motivated by classical Dini conditions on weights (see, e.g., \cite{BloomDeSouza1989}). Notice that, given a weight $\ww=(w_n)_{n=1}^\infty$ and $m\in\NN$,
\begin{equation}\label{eq:Dini}
\sum_{n=1}^m \sum_{k=n}^\infty \frac{w_k}{k}=\sum_{k=1}^\infty \min\{k,m\} \frac{w_k}{k}= \sum_{k=1}^m w_k+m\sum_{k=m+1}^\infty\frac{w_k}{k}.
\end{equation}
Hence, if $\ww$ is nonincreasing and $\uu$ is another nonincreasing weight with primitive sequence $(s_m)_{m=1}^\infty$, the pair $(\uu,\ww)$ is a two-sided Dini pair if and only if
\[
\sup_m \frac{s_m}{m} \sum_{k=1}^\infty \min\{k,m\} \frac{w_k}{k} <\infty.
\]

\begin{proposition}\label{prop:Anso17}
Let $\uu=(u_n)_{n=1}^\infty$ and $\ww=(w_n)_{n=1}^\infty$ be nonincreasing weights with $\lim_n u_n=0$. Set $b=1/u_1$.
\begin{itemize}[leftmargin=*]
\item Assume that 
\begin{equation}\label{eq:minimalcondition}
a:=\sum_{n=1}^\infty \frac{w_n}{n}<\infty.
\end{equation}
Then, $\norm{\cdot}_{\uu,\ww}$ is a quasi-norm on $\BB_{\uu,\ww}$. Furthermore, the following hold.
\begin{enumerate}[label=(\roman*),leftmargin=*,widest=viii]
\item $\norm{\cdot}_{\uu,\ww}$ is equivalent to a $p$-norm for any $0<p<1$.
\item\label{it:B:1} $\norm{T_\psi(\alpha)}_{\uu,\ww}=\norm{\alpha}_{\uu,\ww}$ for all $\alpha\in\co$ and all $\psi\in\Pi(\NN)$.
\item\label{it:B:2} Given $m\in\NN$ and $A\subset \NN$ with $\abs{A}=m$,
\[
\norm{\Ind_A}_{\uu,\ww}=\max\enbrace{\frac{m}{s_m}, \sum_{k=1}^\infty \min\{k,m\} \frac{w_k}{k}}.
\]
\item\label{it:B:3} $\ee_n\in \BB_{\uu,\ww}$ and $\norm{\ee_n}_{\uu,\ww}=\max\{a,b\}$ for all $n\in\NN$.
\item\label{it:B:4} $\restr{\ee_n^*}{\BB_{\uu,\ww}}\in \BB_{\uu,\ww}^*$ and
\[
\frac{1}{\max\{a,b\}} \le \norm{\ee_n^*}_{\BB_{\uu,\ww}\to \FF} \le \frac{1}{b}, \quad n\in\NN.
\]
\item\label{it:B:5} $\BB_{\uu,\ww}$ is complete.
\end{enumerate}

\item In turn, if $(\uu,\ww)$ is a tail Dini pair with constant $K$, then \eqref{eq:minimalcondition} holds, and $\BB_{\uu,\ww}$ has the following additional properties.
\begin{enumerate}[label=(\roman*),leftmargin=*,widest=viii,resume]
\item\label{it:B:6} $\norm{S_A(\alpha)}_{\cG,\ww} \le \norm{\alpha}_{\cG,\ww}+K\norm{\alpha}_{\cM, \uu}$ for all $\alpha\in \co$ and $A\in \Gt(\alpha)$.
\item \label{it:B6b} $(\ee_n)_{n=1}^\infty$ is a quasi-greedy basis of $\BB^0_{\uu,\ww}$ with constant $K+1$.
\item\label{it:B:7} A sequence $\alpha=(a_n)_{n=1}^\infty\in\BB_{\uu,\ww}$ belongs to $\BB_{\uu,\ww}^0$ if and only if
\[
\lim_{A\in\Gt(\alpha)}\norm{\alpha-S_A(\alpha)}_{\cG,\ww} =0
\]
and $\alpha\in m_0(\uu)$.
\end{enumerate}
\end{itemize}
\end{proposition}

\begin{proof}
Set $\epsilon=(2^{1/p-1}-1)/4>0$. Then, $\norm{\cdot}_{\uu,\ww}$ is equivalent to
\[
\norm{\cdot}_\lozenge=\norm{\cdot}_{\cM,\uu} + \norm{\cdot}_{1,\ww} + \epsilon\norm{\cdot}_{\cG,\ww}.
\]
Both maps are positive definite and absolutely homogeneous. Let $\alpha$, $\beta\in\BB_{\uu,\ww}$. Since $\BB_{\uu,\ww}\subset \co$, by Lemma~\ref{lem:Anso1},
\[
\norm{\alpha+\beta}_\lozenge\le (1+4\epsilon)\enpar{ \norm{\alpha}_\lozenge+ \norm{\beta}_\lozenge}=2^{1/p-1}\enpar{ \norm{\alpha}_\lozenge+ \norm{\beta}_\lozenge}.
\]
Hence, $\norm{\cdot}_\lozenge$ is a quasi-norm. Furthermore, by the Aoki--Rolewicz theorem, it is equivalent to a $p$-norm. Consequently, $\norm{\cdot}_{\uu,\ww}$ is a quasi-norm equivalent to a $p$-norm.

Since $m(\uu)$ and $d_1(\ww)$ are symmetric spaces, and the greedy algorithm is invariant under rearrangement, \ref{it:B:1} holds. If $\abs{A}=m$, by \eqref{eq:Gpos} and \eqref{eq:Dini},
\[
\norm{\Ind_A}_{\uu,\ww}
=\max\enbrace{\frac{m}{s_m}, \sum_{n=1}^m w_n, \sum_{n=1}^m \sum_{k=n}^\infty \frac{w_k}{k}}.
\]
Hence, by \eqref{eq:Dini}, \ref{it:B:2} holds. Property~\ref{it:B:3} follows immediately from \ref{it:B:2}. Given $n\in\NN$, since $(\ee_n^*)_{n=1}^\infty$ is a sequence of coordinate functionals for $(\ee_n)_{n=1}^\infty$,
\[
\frac{1}{\max\{a,b\}}\le\norm{\ee_n^*}_{\BB_{\uu,\ww} \to \FF}\le\norm{\ee_n^*}_{m(\uu) \to \FF}=u_1.
\]
Thus, \ref{it:B:4} holds. To prove \ref{it:B:5}, we pick a Cauchy sequence $(\alpha_k)_{k=1}^\infty$. Since $d_1(\ww)$ and $m(\uu)$ are Banach spaces, and $m(\uu)\subset \co$, there is $\alpha\in \co$ such that
\[
\lim_{k\to\infty} \norm{\alpha_k-\alpha}_{\cM,\uu}=\lim_{k\to\infty} \norm{\alpha_k-\alpha}_{1,\ww}=0.
\]
In particular, $\lim_k \alpha_k=\alpha$ uniformly. By Lemma~\ref{lem:Anso2}, for all $k\in\NN$ we have
\[
\norm{\alpha_k-\alpha}_{\cG,\ww}\le \liminf_j \enpar{ \norm{\alpha_k-\alpha_j}_{\cG,\ww} + 2 \norm{\alpha_k-\alpha_j}_{1,\ww}}.
\]
Therefore, $\lim_k \norm{\alpha_k-\alpha}_{\cG,\ww}=0$.

Now suppose that $(\uu, \ww)$ is a tail Dini pair with constant $K$. The definition directly yields
\[
\sum_{n=1}^\infty \frac{w_n}{n}\le w_1+\frac{K}{s_1}<\infty.
\]

Fix $m\in\NN$ and $A\in\Gt(\alpha,m)$. Given $k\in\NN$ and $B\in\Gt(S_A(\alpha),k)$, we estimate
\[
\lambda:=\Ind_B^*(S_A(\alpha))=\Ind_{A\cap B}^*(\alpha)
\]
depending on which of the following exhaustive conditions applies.
\begin{itemize}[leftmargin=*]
\item If $A\subset B$, then $\abs{\lambda}=\abs{\Ind_A^*(\alpha)}\le G_m(\alpha)$.
\item If $\supp(\alpha)\subset A$, then $S_A(\alpha)=\alpha$. Hence, $\abs{\lambda}\le G_k(\alpha)$. 
\item If neither of the preceding cases applies, then $A\subset\supp(\alpha)$ and $A\setminus B\neq\emptyset$. We infer that $B\subset \supp(S_A(\alpha))$. In turn, $\supp(S_A(\alpha))\subset A$. Therefore, $\abs{\lambda}=\abs{\Ind_B^*(\alpha)}\le G_k(\alpha)$.
\end{itemize}

It follows that $G_k(S_A(\alpha))\le G_{\min\{k,m\}}(\alpha)$. Consequently,
\begin{multline*}
\norm{S_A(\alpha)}_{\cG, \ww}\le \sum_{k=1}^m G_k(\alpha) \frac{w_k}{k} +G_m(\alpha) \sum_{k=m+1}^\infty\frac{w_k}{k} \\\
\le \norm{\alpha}_{\cG, \ww}+K \frac{G_m(\alpha)}{s_m} \le \norm{\alpha}_{\cG, \ww}+K\norm{\alpha}_{\cM, \uu}.
\end{multline*}

Hence, \ref{it:B:6} holds. Combining \ref{it:B:6} with the fact that the unit vector system $\Et=(e_n)_{n=1}^\infty$ is $1$-unconditional relative to $m(\uu)$ and $d_1(\ww)$ gives \ref{it:B6b}. Hence, by Theorem~\ref{thm:AABWQG}, a sequence $\alpha\in\co$ belongs to $\BB^0_{\uu,\ww}$ if and only if
\[
\lim_{A\in\Gt(\alpha)}\norm{\alpha-S_A(\alpha)}_{\uu,\ww} =0.
\]
Since the linear span of $\Et$ is dense in $\co \cap d_1(\ww)$, \ref{it:B:7} holds.
\end{proof}

 Our proof of Theorem~\ref{thm:BeraA} relies on the existence of two-sided Dini pairs satisfying special properties.

\begin{definition}
Let $\uu=(u_n)_{n=1}^\infty$ and $\ww=(w_n)_{n=1}^\infty$ be nonincreasing weights, and let $(s_m)_{m=1}^\infty$ be the primitive sequence of $\uu$. We say that $(\uu,\ww)$ is a \emph{fluttery pair} if it is a two-sided Dini pair,
\[
\sum_{n=1}^\infty \frac{w_n s_n }{n}=\infty,
\]
and
\[
\sum_{n=1}^\infty u_n w_n<\infty.
\]
\end{definition}

\begin{remark}
If the weight $\ww=(w_n)_{n=1}^\infty$ satisfies the tail Dini condition
\[
\sum_{n=m}^\infty\frac{w_n}{n} \approx w_m, \quad m\in\NN,
\]
then, by the Abel summation formula,
\[
\sum_{n=1}^\infty \frac{w_n s_n }{n}= \sum_{m=1}^\infty u_m \sum_{n=m}^\infty \frac{w_n}{n} \sim \sum_{m=1}^\infty u_m w_m.
\]
Hence, since Dini conditions are widely regarded as regularity conditions on weights, the existence of $\uu$ such that $(\uu,\ww)$ is fluttery points to a lack of regularity of $\ww$.
\end{remark}

\begin{lemma}\label{lem:unull}
Let $\uu=(u_n)_{n=1}^\infty$ and $\ww=(w_n)_{n=1}^\infty$ be weights. If $(\uu,\ww)$ is a fluttery pair, then $\lim_n u_n=0$.
\end{lemma}

\begin{proof}
Assume, towards a contradiction, that $\lim_n u_n>0$. Then $u_n\approx 1$ for $n\in\NN$, and the primitive sequence $(s_m)_{m=1}^\infty$ of $\uu$ satisfies $s_m\approx m$ for $m\in\NN$. Consequently, 
\[
\sum_{n=1}^\infty \frac{w_n s_n }{n}\sim \sum_{n=1}^\infty w_n \sim \sum_{n=1}^\infty u_n w_n.
\]
This absurdity proves that $\lim_n u_n=0$.
\end{proof}

\begin{example}\label{ex:flut}
Define $\uu=(u_n)_{n=1}^\infty$ and $\ww=(w_n)_{n=1}^\infty$ by
\[
u_n=\frac{1}{n},\quad w_n=\frac{1}{\log^2 (1+n)}.
\]
Since $s_m=\sum_{n=1}^m u_n\approx \log(1+m)$ for $m\in\NN$,
\[
\sum_{n=1}^m w_n \approx \int_2^{m+2} \frac{ds}{\log^2(s)} \approx \frac{m+2}{\log^2(2+m)}
\lesssim \frac{m}{\log(1+m)},
\quad m\in\NN,
\]
and
\[
\sum_{n=m}^\infty \frac{w_n}{n} \approx \int_{m+1}^\infty \frac{ds}{s\log^2(s)} = \frac{1}{\log(m+1)},
\quad m\in\NN,
\]
$(\uu, \ww)$ is a fluttery pair.
\end{example}

\begin{proof}[Completion of the proof of Theorem~\ref{thm:BeraA}]
Choose a fluttery pair $(\uu,\ww)$ with $\uu$ decreasing, whose existence is guaranteed by Example~\ref{ex:flut}. Define $(s_m)_{m=1}^\infty$ as the primitive sequence of $\uu=(u_n)_{n=1}^\infty$. Let $\YB$ be the canonical basis of $m_0(\uu)$. By Lemma~\ref{lem:unull}  and Proposition~\ref{prop:Anso17}, $\BB_{\uu,\ww}$ is a quasi-Banach space and the unit vector system $\XB$ is an M-bounded Markushevich basis of $\BB^0_{\uu,\ww}$. Moreover, $\XB$ is quasi-greedy and symmetric, and
\begin{multline*}
\udf[\XB](m) =\sup_{\abs{A}\le m} \norm{\Ind_A}_{\uu,\ww}=\max\enbrace{\frac{m}{s_m}, \sum_{k=1}^\infty \min\{k,m\} \frac{w_k}{k}} \\
\approx \frac{m}{s_m}=\udf[\YB](m).
\end{multline*}
By Proposition~\ref{prop:SubSymBid}, $\YB$ is bidemocratic. Consequently, $\XB$ is bidemocratic by Lemma~\ref{lem:alt_norm_bidem}.

It remains to prove that no rearrangement turns $\XB$ into a Schauder basis. By Lemma~\ref{lem:ScSymUnc}, it suffices to show that $\XB$ is not unconditional. Set $\uu=(u_n)_{n=1}^\infty$. Define for each $j\in\NN$
\[
\alpha_j=\sum_{n=1}^j u_n\, \ee_n, \quad \beta_j=\sum_{n=1}^j (-1)^n u_n\, \ee_n.
\]
We have $G_m(\alpha_j)=s_{\min\{j,m\}}$ and, by the Leibniz criterion for alternating series, $G_m(\beta_j)\le u_1$ for all $j$, $m\in\NN$. Hence,
\[
\norm{\beta_j}_{\cG,\ww}\le u_1 \sum_{n=1}^\infty \frac{w_n}{n}
\]
for all $j\in\NN$. Moreover,
\[
\norm{\beta_j}_{\cM,\uu}=1, \quad \norm{\beta_j}_{1,\ww}=\sum_{n=1}^j u_n w_n.
\]
Thus, $(\beta_j)_{j=1}^\infty$ is bounded in $\BB_{\uu,\ww}$. In contrast, by \eqref{eq:Gpos},
\[
\norm{\alpha_j}_{\uu,\ww} \ge \norm{\alpha_j}_{\cG,\ww}
\ge \sum_{m=1}^j \frac{w_m}{m} G_m(\alpha_j)=\sum_{m=1}^j \frac{w_m s_m}{m}.
\]
Consequently, $(\alpha_j)_{j=1}^\infty$ is unbounded in $\BB_{\uu,\ww}$, and we are done.
\end{proof}

We have proven that bidemocracy and symmetry are not strong enough conditions to imply unconditionality in the non-locally convex setting. We conclude by showing that they do entail it in the Banach envelope.
\begin{lemma}\label{lem:SymEnv}
Let $\XB$ be a fundamental minimal system of a quasi-Banach space $\XX$. Assume that $\XB$ is symmetric. Then, the envelope fundamental minimal system $\XB_c$ is symmetric.
\end{lemma}

\begin{proof}
Let $\psi\in\Pi(\NN)$. Since $T_\psi[\XB]$ is a bounded operator from $\XX$ to $\XX$, $T_\psi[\XB_c]$ is a bounded operator from $\XX_c$ to $\XX_c$.
\end{proof}

\begin{proposition}\label{prop:SymBdEnv}
Let $\XB$ be a fundamental minimal system of a quasi-Banach space $\XX$. Assume that $\XB$ is symmetric and bidemocratic. Then, the envelope fundamental minimal system $\XB_c$ is unconditional.
\end{proposition}

\begin{proof}
This follows by combining Lemma~\ref{lem:bidem_envelope}, Lemma~\ref{lem:SymEnv} and Theorem~\ref{thm:dem-spread-equivalences}.
\end{proof}

Proposition~\ref{prop:SymBdEnv} highlights a relevant difference between the quasi-greedy non-Schauder system defined in \cite{AABB2026}*{Section 4} and that provided by Theorem~\ref{thm:BeraA}. In the former example, quasi-greediness and even totality are lost upon passage to the Banach envelope, whereas in the latter example passing to the Banach envelope improves the greedy-like properties of the system.

\section*{Statements and Declarations}
\subsection*{Conflict of Interest}
The authors declare that they have no conflict of interest.
\subsection*{Data Availability}
Since no datasets were generated or analyzed during the current study, data sharing does not apply to this article.
\subsection*{Generative AI use}
This paper was conceived, developed, and written without using any AI tools. The authors are solely responsible for the mathematical results and their proofs.
\bibliography{references}
\end{document}